\documentclass[10pt]{amsart} 
\usepackage{amsmath} 
\usepackage{amsxtra}
\usepackage{amscd} 
\usepackage{amsthm} 
\usepackage{amsfonts}
\usepackage{amssymb} 
\usepackage{eucal}

\newtheorem{cor}[subsection]{Corollary}
\newtheorem{lem}[subsection]{Lemma}
\newtheorem{prop}[subsection]{Proposition}

\newtheorem{conj}[subsection]{Conjecture}
\newtheorem{thm}[subsection]{Theorem}

\theoremstyle{definition}

\theoremstyle{remark}

\newcommand{\thmref}[1]{Theorem~\ref{#1}}
\newcommand{\secref}[1]{Sect.~\ref{#1}}
\newcommand{\lemref}[1]{Lemma~\ref{#1}}
\newcommand{\propref}[1]{Proposition~\ref{#1}}
\newcommand{\corref}[1]{Corollary~\ref{#1}}
\newcommand{\conjref}[1]{Conjecture~\ref{#1}}

\newcommand{\nc}{\newcommand}
\nc{\renc}{\renewcommand}
\nc{\ssec}{\subsection}
\nc{\sssec}{\subsubsection}
\nc{\on}{\operatorname}

\nc\ol{\overline}
\nc\wt{\widetilde}
\nc\tboxtimes{\wt{\boxtimes}}
\nc{\alp}{\alpha}
\nc{\hl}{\overset{\leftarrow}h}
\nc{\hr}{\overset{\rightarrow}h}

\nc{\ZZ}{{\mathbb Z}}
\nc{\NN}{{\mathbb N}}
\nc{\PP}{{\mathbb P}}
\nc{\FF}{{\mathbb F}}
\nc{\OO}{{\mathbb O}}
\renc{\SS}{{\mathbb S}}
\nc{\DD}{{\mathbb D}}
\nc{\GG}{{\mathbb G}}

\nc{\Fq}{{\mathbb F}_q}
\nc{\Fqb}{\ol{{\mathbb F}_q}}
\nc{\Ql}{\ol{{\mathbb Q}_\ell}}
\nc{\id}{\text{id}}
\nc\X{\mathcal X}

\nc{\Hom}{\on{Hom}}
\nc{\Lie}{\on{Lie}}
\nc{\Loc}{\on{Loc}}
\nc{\Pic}{\on{Pic}}
\nc{\Bun}{\on{Bun}}
\nc{\IC}{\on{IC}}
\nc{\Aut}{\on{Aut}}
\nc{\rk}{\on{rk}}
\nc{\Sh}{\on{Sh}}
\nc{\Perv}{\on{Perv}}
\nc{\pos}{{\on{pos}}}
\nc{\Conv}{\on{Conv}}
\nc{\Sph}{\on{Sph}}
\nc{\Sym}{\on{Sym}}
\nc{\sym}{\on{sym}}
\nc{\BunBb}{\overline{\Bun}_B}
\nc{\Buno}{\overset{o}{\Bun}}
\nc{\BunPb}{{\overline{\Bun}_P}}
\nc{\BunBM}{\overline{\Bun}_{B(M)}}
\nc{\BunPbw}{{\widetilde{\Bun}_P}}
\nc{\BunBP}{\widetilde{\Bun}_{B,P}}
\nc{\GUb}{\overline{G/U}}
\nc{\GUPb}{\overline{G/U(P)}}
\nc{\Dt}{\widetilde{\on{D}}}
\nc{\Pt}{\widetilde{\on{P}}}

\nc{\Hhom}{\underline{\on{Hom}}}
\nc\syminfty{\on{Sym}^{\infty}}
\nc\lal{\ol{\lambda}}
\nc\xl{\ol{x}}
\nc\thl{\ol{\theta}}
\nc\nul{\ol{\nu}}
\nc\mul{\ol{\mu}}
\nc\Sum\Sigma

\nc{\oX}{\overset{o}{X}{}}

\nc{\M}{{\mathcal M}}
\nc{\N}{{\mathcal N}}
\nc{\F}{{\mathcal F}}
\nc{\D}{{\mathcal D}}
\nc{\Q}{{\mathcal Q}}
\nc{\Y}{{\mathcal Y}}
\nc{\G}{{\mathcal G}}
\nc{\CG}{{\mathcal G}}
\nc{\CC}{{\mathcal C}}
\nc{\E}{{\mathcal E}}
\nc{\CalC}{{\mathcal C}}
\nc\Dh{\widehat{\D}}
\renewcommand{\O}{{\mathcal O}}
\nc{\C}{{\mathcal C}}
\nc{\K}{{\mathcal K}}
\renewcommand{\H}{{\mathcal H}}
\renewcommand{\S}{{\mathcal S}}
\nc{\T}{{\mathcal T}}
\nc{\V}{{\mathcal V}}
\renc{\P}{{\mathcal P}}
\nc{\A}{{\mathcal A}}
\nc{\B}{{\mathcal B}}
\nc{\U}{{\mathcal U}}
\renewcommand{\L}{{\mathcal L}}
\nc{\Gr}{\on{Gr}}

\nc{\fA}{{\mathfrak A}}
\nc{\fP}{{\mathfrak P}}
\nc{\frn}{{\check{\mathfrak u}(P)}}
\nc{\p}{\mathfrak p}
\nc{\q}{\mathfrak q}
\nc\f{{\mathfrak f}}

\nc{\s}{{\mathfrak s}}
\nc\w{\text{w}}

\nc\mathi\iota
\nc\Spec{\on{Spec}}
\nc\Mod{\on{Mod}}
\nc{\tw}{\widetilde{\mathfrak t}}
\nc{\pw}{\widetilde{\mathfrak p}}
\nc{\qw}{\widetilde{\mathfrak q}}
\nc{\jw}{\widetilde j}

\nc{\grb}{\overline{\Gr}}
\nc{\I}{\mathcal I}

\nc{\lambdach}{\check\lambda}
\nc{\Lambdach}{\check\Lambda}
\nc{\much}{\check\mu}
\nc{\omegach}{\check\omega}
\nc{\nuch}{\check\nu}
\nc{\etach}{\check\eta}
\nc{\alphach}{\check\alpha}
\nc{\betach}{\check\beta}
\nc{\rhoch}{\check\rho}
\nc{\wh}{\widehat}

\nc{\BZ}{{\mathbb Z}}
\nc{\BF}{{\mathbb F}}
\nc{\bF}{{\mathbf F}}
\nc{\bO}{{\mathbf O}}

\nc{\Hb}{\overline{\H}}

\nc{\Qb}{\ol{\Q}}
\nc{\db}{\ol{d}}
\nc{\yb}{\ol{y}}

\nc{\Dwk}{\on{D}^W(\Qb_k)}
\nc{\Dwkex}{\on{D}^W(\Qb_{k,ex})}
\nc{\Dwkexp}{\on{D}^W(\Qb_{k+1,ex})}
\nc{\Pwk}{\on{P}^W(\Qb_k)}
\nc{\Pwkex}{\on{P}^W(\Qb_{k,ex})}
\nc{\Pwkexp}{\on{P}^W(\Qb_{k+1,ex})}

\nc{\Rep}{\on{Rep}}
\nc{\cG}{{\check{G}}}
\nc{\CT}{\mathcal T}
\nc{\CR}{\mathcal R}
\nc{\Ind}{\on{Ind}}
\nc{\sF}{\mathsf F}
\nc{\sG}{\mathsf G}
\nc{\CK}{\mathcal K}
\nc{\CO}{\mathcal O}
\nc{\oDelta}{\overset{o}\Delta{}}
\nc{\oGr}{\overset{o}{\Gr}{}}
\nc{\oR}{\overset{o}\CR{}}
\nc{\fC}{\mathfrak C}
\nc{\shriektimes}{\overset{!}\otimes}

\title{On de Jong's conjecture}

\author{D. Gaitsgory}

\address{Dept. of Math., The University of Chicago, 
5734 South University Ave., Chicago, IL, 60637}

\email{gaitsgde@math.uchicago.edu}

\begin{document}

\maketitle

\section{Introduction}

\ssec{}

The purpose of this note is to indicate the proof of a partial case of 
de Jong's conjecture,
proposed in \cite{DJ}. Let us recall its formulation, combining Conjecture 2.3
and Theorem 2.17 of {\it loc. cit.}:

Let $X$ be a smooth projective curve over a finite field $\BF_q$,
and let $\rho$ be a continuous representation
$$\pi_1(X)\to GL_n(\bF),$$
where $\bF=\BF_l((t))$ with $\BF_l$ being another finite 
field of order coprime to $q$.
Assume that $\rho|_{\pi_1(\overline{X})}$ is absolutely irreducible
(here, as usual, $\pi_1(\overline{X})\subset \pi_1(X)$ denotes the
geometric fundamental group).

\begin{conj}   \label{De Jong}
Under the above circumstances, $\rho(\pi_1(\overline{X}))$ is
finite.
\end{conj}

In the same paper, de Jong showed that \conjref{De Jong} follows
from a version of the Langlands conjecture with $\bF$-coefficients.
(The corresponding version of the Langlands conjecture with $\Ql$-coefficients 
is now a theorem of Lafforgue.) 

Namely, consider the double quotient $GL_n(K)\backslash GL_n({\mathbb A})/GL_n({\mathbb O})$,
where $K$ is the global field corresponding to $X$, ${\mathbb A}$ is the 
ring of adeles, and ${\mathbb O}$ is the subring of integral adeles. (Note 
that the above double quotient 
identifies with the set of isomorphism classes of rank $n$ vector bundles
on $X$, denoted $\Bun_n(\BF_q)$.) For each place $x\in |X|$ one introduces the Hecke 
operators $T^i_x$, $i=1,...,n$ acting on the space of $\bF$-valued functions on 
$\Bun_n(\BF_q)$.

Given a representation $\rho: \pi_1(X)\to GL_n(\bF)$,
we say that a function $f:\Bun_n({\mathbb F}_q)\to \bF$
is a Hecke eigen-form with eigenvalues corresponding to $\rho$, if for all $x$ and $i$
$$T^i_x(f)=\lambda^i_x\cdot f,$$
with $\lambda^i_x=\on{Tr}(\Lambda^i(\rho(Fr_x)))$, where $Fr_x\in \pi_1(X)$ is the 
Frobenius element corresponding to $x$ (defined up to conjugacy), and $\Lambda^i$
designates the $i$-th exterior power of the representation $\rho$.

In addition, in the space of all $\bF$-valued functions on 
$\Bun_n({\mathbb F}_q)$ one singles out a subspace,
stable under the operators $T_x^i$, of cuspidal functions. As in the $n=2$ case,
one shows that the subspace of cuspidal functions with a given central character
is finite-dimensional.

Here is the relevant version of the Langlands conjecture:

\begin{conj}  \label{Langlands}
Let $\rho$ be as in \conjref{De Jong}. Then there exists a (non-zero)
cuspidal Hecke-eigenform $f_\rho$ with eigen-values corresponding to
$\rho$.
\end{conj}

By arguments of Sect. 4 of \cite{DJ}, one shows that \conjref{Langlands} implies
\conjref{De Jong}. In this paper we will be concerned with the proof of 
\conjref{Langlands}. 

\ssec{}

Unfortunately, the proof of \conjref{Langlands} is not complete.
Namely, we will have to rely on two pieces of mathematics
that do not exist in the published literature.

One is the theory of \'etale sheaves with $\bF$ coefficients, 
which should be parallel to the theory of $\bF'$-adic sheaves, 
where $\bF'$ is a local field of characteristic $0$.

However, we do not expect that the
construction of this theory is in any way different from its $\bF'$-counterpart.
Namely, we first consider \'etale sheaves of ${\mathbb F}_l[t]/t^i$-modules,
and define the category of $\bO$-sheaves as the appropriate 2-projective 
limit, where $\bO$ is the local ring $\BF_l[[t]]$.
We define the category of $\bF$-sheaves by inverting $t$, i.e., by quotienting out 
the category of $\bO$-sheaves by $t$-torsion ones.

Thus, from now on we will assume the existence of such a theory, with formal
properties analogous to that of $\bF'$-sheaves. In particular, we assume the
existence of the category $D(\Y)$ of bounded complexes with constructible
cohomology corresponding to an algebraic variety $\Y$ over ${\mathbb F}_q$, stable under the
6-functors. From this category one produces the abelian category of
perverse sheaves, denoted $P(\Y)$. 

Finally, we must have the sheaf-function correspondence. In other words,
given an object $\K\in D(\Y)$, by taking traces of the Frobenius elements, 
we obtain a $\bF$-valued function on $\Y({\mathbb F}_q)$, and this operation is 
compatible with the operations of taking inverse image, direct image with
compact supports and tensor product.

\ssec{}

That said, our goal will be to prove a geometric version of \conjref{Langlands}
(we refer the reader to \secref{formulation of theorems} for the precise
formulation of the latter).

Namely, to a representation $\rho$, as in \conjref{De Jong}, we would like 
to associate an object $\S_\rho\in D(\Bun_n)$, where $\Bun_n$ denotes the
moduli stack of rank $n$ bundles on $X$, such that $\S_\rho$ is cuspidal
and satisfies the Hecke eigen-condition with respect to $\rho$.
If such $\S_\rho$ exists, then by applying the sheaf-function correspondence,
we obtain a function $f_\rho$ on $\Bun_n({\mathbb F}_q)$, which is a cuspidal
Hecke eigen-form with eigenvalues corresponding to $\rho$.

Our main result is \thmref{main} that says that if $l>2n$, and $\rho$ is
as in \conjref{De Jong}, then the object $\S_\rho\in D(\Bun_n)$ with the
required properties exists.

In fact, we state and indicate the proof of a stronger result, namely,
\thmref{strong form}, which asserts the existence of $\S_\rho$ for
any $l\neq 2$. 

To summarize, this paper proves \conjref{De Jong} for $l\neq 2$ modulo
the theory of $\bF$-sheaves, and one more unpublished result,
discussed below.

\ssec{}

Even in order to formulate the geometric analog of \conjref{Langlands}, i.e.,
\conjref{geometric langlands}, one needs to rely on the realization of the
category of representations of the Langlands dual group via spherical perverse
sheaves on the affine Grassmannian.

This result was first announced by V.~Ginzburg in 1995's for 
perverse sheaves with coefficients in a field of characteristic zero.
Recently, in \cite{MV}, I.~Mirkovi\'c and K.~Vilonen have established 
this result in a far greater generality 
(cf. Theorem 12.1 in {\it loc. cit}). 

Namely, they work over a {\it ground field} of complex numbers,
and consider sheaves in the analytic topology. They show that
the category of spherical perverse sheaves with $A$-coefficients,
where $A$ is an arbitrary Noetherian commutative ring of finite
cohomological dimension, is equivalent to the category of 
representations of the group-scheme $\check G_A$ on finitely
generated $A$-modules, where $\check G_A$ is the split reductive
group, whose root datum is dual to that of $G$.

\medskip

For the purposes of this paper, however, we need an extension
of the result of \cite{MV} to the case of an arbitrary ground field 
(in practice taken to be $\BF_q$), and sheaves with $\bF$-
and $\bO$-coefficients. From examining \cite{MV} it appears that 
the proof presented in {\it loc. cit.} carries over to this case. 
Therefore, we state the corresponding result as \thmref{satake}.

\ssec{}

The geometric Langlands conjecture with coefficients of
characteristic $0$ has been proved in \cite{FGV} and \cite{Ga}.
The proof of \thmref{main} follows verbatim the approach 
of {\it loc.cit.}, with one exception:

This exception is the discussion related to the notion of
symmetric power of a local system (cf. \secref{sect powers}),
and this is the only original piece of work done in this 
paper.

\medskip

Let us now describe the contents of the paper:

\medskip

In \secref{rev satake} we recall the definition of affine
Grassmannians and review the relevant versions of the 
geometric Satake equivalence, i.e., the realization of the 
Langlands dual group via spherical perverse sheaves, following
\cite{MV}. We also introduce the Hecke stack, Hecke
functors and the notion of Hecke eigen-sheaf.

\medskip

Starting from \secref{intr conj} we restrict ourselves to the case
$G=GL_n$. In the beginning of \secref{intr conj} we formulate
several versions of the geometric Langlands conjecture for $GL_n$,
and state the main result, \thmref{main}, which amounts to 
proving one of the forms of the above conjecture when 
$\on{char}(\bF)>2n$. We also formulate a stronger result, 
\thmref{strong form}, which claims the ultimate form of the
geometric Langlands conjecture when $\bF\neq 2$.

In the rest of \secref{intr conj} we explain, following closely the
exposition in \cite{FGV}, how \thmref{main} can be reduced to a certain
vanishing statement, \thmref{vanishing}.

\medskip

In \secref{sect vanishing} we formulate and prove \thmref{vanishing},
following \cite{Ga}. The proof essentially amounts to introducing a
certain quotient triangulated category $\wt{D}(\Bun_n)$ of
$D(\Bun_n)$ and showing that the averaging functor 
$\on{Av}^d_E:D(\Bun_n)\to D(\Bun_n)$ (whose vanishing for large $d$
we are trying to prove) is well-defined and {\it exact} on the
above quotient.

The proof of the latter fact essentially consists of two steps:

\smallskip

\noindent (1) Showing that the elementary functor 
$\on{Av}^1_E:\wt{D}(\Bun_n)\to \wt{D}(\Bun_n)$ is exact.

\smallskip

\noindent (2) Showing that the exactness of $\on{Av}^1_E$
implies the exactness of $\on{Av}^d_E$ for any $d$.

\smallskip

Step (1) requires no modification compared to the case 
of characteristic $0$ coefficients if $\on{char}(\bF)>2n$
considered in \cite{Ga}, and a minor modification if we only 
assume $\on{char}(\bF)\neq 2$ (the latter case is treated 
in the Appendix).

Step (2) amounts to \propref{functor on quotient}(2) and
this is the only place in the paper that requires some
substantial work. Essentially, $\on{Av}^d_E$
is the $d$-th symmetric power (along $X$) of $\on{Av}^1_E$,
and our goal is to express one through another. This is 
achieved by introducing a somewhat non-standard notion of
external exterior power of a local system on a curve.

\medskip

In \secref{sect powers} we first recall some notions from
linear algebra, namely, the two versions of symmetric and
exterior powers of a vector space and the corresponding
Koszul complexes, when working over a field of positive
characteristic. 

We then review some basic properties of
the construction of the (external) symmetric power of 
a local system on a curve, in particular, its behavior with
respect to the perverse t-structure. 

And finally, we introduce
two versions of an external exterior power of a local system,
and construct the corresponding external Koszul complex,
which is used in the proof of \propref{functor on quotient}(2) mentioned above.

In Appendix A we indicate how, by refining
some arguments of \cite{Ga}, one can relax the condition
that $\on{char}(\bF)>2n$ and treat the case of any $\bF$
of characteristic different from $2$.

Finally, in Appendix B we prove \thmref{extended satake}, which
is a version of the geometric Satake equivalence over
a symmetric power of the curve $X$.

\ssec{Acknowledgments}

I would like to thank V.~Drinfeld for many helpful discussions.

This work was supported by a long-term fellowship at the 
Clay Mathematics Institute, and by a grant from DARPA  
via NSF, DMS 0105256.

\section{Review of the geometric Satake equivalence}  \label{rev satake}

From now on, we will work over an arbitrary algebraically closed 
ground field $k$ of characteristic prime to $l$, 
and $X$ will be a smooth projective curve $/k$.

In this section $G$ will be an arbitrary reductive group over $k$.
By $\check G$ we will denote its Langlands dual group, which we 
think of as a smooth group-scheme defined over $\BZ$. We will
denote by $\check G_\bO$ the corresponding group-scheme
over $\bO$ and by $\check G_\bF$ the corresponding reductive 
group over $\bF$. We will denote by $\on{Rep}(\check G_\bF)$ the
category of rational representations of $\check G_\bF$ on
finite-dimensional $\bF$-vector spaces.

\ssec{}

Let $x\in X$ be a point. If $G$ is an algebraic group,
let $\Gr_{G,x}$ denote the affine Grassmannian of $G$ on $X$ at $x$. In
other words, $\Gr_{G,x}$ is the ind-scheme classifying the data of a
$G$-bundle $\P_G$ on $X$ with a trivialization $\beta$ on $X-x$, i.e., 
$\P_G|_{X-x}\simeq \P^0_G|_{X-x}$, where $\P^0_G$ denotes the trivial $G$-bundle.

According to a theorem of Beauville and Laszlo, the data of $(\P_G,\beta)$ 
is equivalent to one, where instead of $X$ we use the formal disc $\D_x$ around $x$,
and instead of $X-x$ the formal punctured disc $\D^*_x$, cf. \cite{Ga1}.

Let $G(\K_x)$ (resp., $G(\O_x)$) be the group ind-scheme (resp., group-scheme)
classifying maps $\D^*_x\to G$ (resp., $\D_x\to G$). The description of $\Gr_{G,x}$
via $\D_x$ implies that the group ind-scheme $G(\K_x)$ acts on it
by changing the data of $\beta$. The action of $G(\O_x)$ has the property that
every finite-dimensional closed subscheme of $\Gr_{G,x}$ is contained in
another finite-dimensional closed subscheme of $\Gr_{G,x}$,
stable under the action of $G(\O_x)$. Therefore, the category $P^{G(\O_x)}(\Gr_{G,x})$
of $G(\O_x)$-equivariant perverse sheaves (with $\bF$-coefficients)
on $\Gr_{G,x}$ makes sense.

The basic fact is that $P^{G(\O_x)}(\Gr_{G,x})$ carries a natural structure of
monoidal category. We refer to \cite{Ga1}, where this is discussed in detail.

We will need the following result, which is a generalization of Theorem 12.1 of
\cite{MV} to the case of an arbitrary ground field:

\begin{thm}  \label{satake}
The monoidal structure on $P^{G(\O_x)}(\Gr_{G,x})$ admits a natural symmetric
commutativity constraint. The resulting tensor category is equivalent to 
the category $\on{Rep}(\check G_\bF)$.
\end{thm}

Let $\on{Aut}(\D_x)$ be the group-scheme of automorphisms of the formal disc
$\D_x$. By functoriality, we have a natural action of $\on{Aut}(\D_x)$ on
$\Gr_{G,x}$. A part of \thmref{satake} is the following statement
(cf. Proposition 2.2 of \cite{MV}):

\begin{cor} \label{virasoro}
Every object of $P^{G(\O_x)}(\Gr_{G,x})$ is $\on{Aut}(\D_x)$-equivariant.
\end{cor}

\medskip

\thmref{satake} remains valid in the context of $\bO$-coefficients.
Namely, let us denote by $P^{G(\O_x)}(\Gr_{G,x})_\bO$ the 
category of $\bO$-valued $G(\O_x)$-equivariant perverse sheaves
on $\Gr_{G,x}$. \footnote{Since $\bO$ is a ring and not a field, Serre
duality on the derived category of $\bO$-modules does not preserve
the t-structure. Therefore, the corresponding derived category of
$\bO$-sheaves (as in the ${\mathbb Z}_l$ case) possesses two 
natural perverse t-structures, interchanged by Verdier duality.
Here we will use the "usual" one, for which the category of perverse 
sheaves in Noetherian.}
Then $P^{G(\O_x)}(\Gr_{G,x})_\bO$ is also
a tensor category, equivalent to the category of representations
of the group-scheme $\check G_\bO$ on finitely generated
$\bO$-modules.

\ssec{}

We will now discuss a version of \thmref{satake}, where instead of one point $x$
we have several points moving along the curve $X$. For a positive integer $d$, 
let $\Gr^d_G$ be the symmetrized version of the Beilinson-Drinfeld Grassmannian living over
$X^{(d)}$. Namely, a point of $\Gr^d_G$ is a triple $(D,\P_G,\beta)$, where $D$ is an effective
divisor of degree $d$ on $X$, i.e., a point of $X^{(d)}$, $\P_G$ is a principal
$G$-bundle on $X$, and $\beta$ is a trivialization of $\P_G$ off the support of $D$.

The formal disc version of $\Gr^d_G$ can be spelled out as follows. 
Let $S$ be a (test) scheme and let $D_S$ be an $S$-point of $X^{(d)}$. Let
$\Gamma$ be the incidence divisor in $X\times X^{(d)}$, and let $\Gamma_S$
be its pull-back to $X\times S$. Let $k\cdot \Gamma_S$ denote the
$k$-th infinitesimal neighborhood of $\Gamma_S$ and let $\wh{\Gamma}_S$ be the 
completion of $X\times S$ along $\Gamma_S$. Note that it makes sense to speak
about principal $G$-bundles on $\wh{\Gamma}_S$, and of isomorphisms of two
such bundles on $\wh{\Gamma}_S^*:=\wh{\Gamma}_S-\Gamma_S$.

Thus, an $S$-point of $\Gr^d_G$ is a triple $(D_S,\P_G,\beta)$, where
$D_S$ is an $S$-point of $X^{(d)}$, $\P_G$ is a principal $G$-bundle on
$\wh{\Gamma}_S$, and $\beta$ is an isomorphism $\P_G\simeq \P^0_G$ on
$\wh{\Gamma}_S^*$.

\medskip

For a partition $\db:d=d_1+...+d_m$, let $sum_{\db}$ denote the natural morphism
$X^{(d_1)}\times...\times X^{(d_m)}\to X^{(d)}$. Let
$X^{(\db)}_{disj}$ denote the open subset of
$X^{(d_1)}\times...\times X^{(d_m)}$ corresponding to the condition that 
all the divisors $D_i\in X^{(d_i)}$, $i=1,...,m$ have pairwise disjoint
supports. We have:

\begin{equation}  \label{factorization}
X^{(\db)}_{disj}\underset{X^{(d)}}\times \Gr^d_G\simeq  
X^{(\db)}_{disj}\underset{X^{(d_1)}\times...\times X^{(d_m)}}\times
(\Gr^{d_1}_G\times...\times \Gr^{d_m}_G).
\end{equation}

\medskip

For an integer $k$, let $\CG^d_k$ denote the group-scheme over $X^{(d)}$,
whose $S$-points over a given $D_S\in \on{Hom}(S,X^{(d)})$ is the group
of maps $k\cdot \Gamma_S\to G$. Let $\CG^d$ be the group-scheme
$\underset{\longleftarrow}{lim}\, \CG^d_k$. We have a natural action of
$\CG^d$ on $\Gr^d_G$.
Note that the fiber of $\CG^d$ at $D=\Sigma\, d_i\cdot x_i$ with $x_i$'s distinct
is $\Pi\, G(\O_{x_i})$. These isomorphisms are easily seen to be compatible with
the factorization isomorphism of \eqref{factorization}.

We let $P^{\CG^d}(\Gr^d_G)$ denote the category of
$\CG^d$-equivariant perverse sheaves on $\Gr^d_G$. 
We claim that for $d=d_1+d_2$ there exists a natural convolution functor
\begin{equation} \label{conv of sheaves}
\star:P^{\CG^{d_1}}(\Gr^{d_1}_G)\times P^{\CG^{d_2}}(\Gr^{d_2}_G)\to P^{\CG^d}(\Gr^d_G).
\end{equation}

Indeed, consider the ind-scheme $\on{Conv}_G^{d_1,d_2}$ that classifies the data of
$$(D_1,D_2,\P^1_G,\beta^1,\P^G,\beta),$$ where $D_i\in X^{(d_i)}$, $\P^1_G,\P_G$ are
principal $G$-bundles on $X$, $\beta_1:\P^1_G|_{X-D_1}\simeq \P_G|_{X-D_1}$,
$\beta:\P_G|_{X-D_2}\simeq \P^0_G|_{X-D_2}$. We can view $\on{Conv}_G^{d_1,d_2}$ 
as a fibration over $\Gr^{d_2}_G$ with the typical fiber $\Gr^{d_1}_G$. Thus, for
any perverse sheaf $\T_2$ on $\Gr^{d_2}_G$ and a $\CG^{d_1}$-equivariant perverse sheaf
$\T_1$ on $\Gr^{d_1}_G$ we can associate their twisted external product
$\T_1\tboxtimes \T_2\in P(\on{Conv}_G^{d_1,d_2})$.  
We also have a natural projection $p:\on{Conv}_G^{d_1,d_2}\to 
\Gr^d_G$ that sends a data $(D_1,D_2,\P^1_G,\beta^1,\P^G,\beta)$ as above to
$(D_1+D_2,\P^1_G,\beta\circ\beta^1)$. We set
$$\T_1\star \T_2:=p_!(\T_1\tboxtimes \T_2).$$

The fact that the resulting object of the derived category is a perverse sheaf
follows from the semi-smallness of convolution, cf. \cite{MV}, Sect. 4.1. The fact that 
this perverse sheaf if $\CG^d$-equivariant is evident, since all the objects involved in the
construction, when viewed over $X^{(d)}$, carry a natural $\CG^d$-action.

\ssec{}

Note that if $\CC$ is an $\bF$-linear abelian category, it makes sense to 
speak about objects of $\CC$, 
endowed with an action of the algebraic group $\check G_\bF$. 
We will take $\CC$ to be the category of $\bF$-perverse sheaves on various schemes. 

We introduce the category $P^{\check G,d}$ to consist of perverse sheaves
$\K$ on $X^{(d)}$, endowed with the following structure:
For any ordered partition $\db:d=d_1+...+d_m$, the
pull-back
$$\K^{\db}:=sum_{\db}^*(\K)|_{X^{(\db)}_{disj}}$$
carries an action of $(\check G_\bF)^{\times m}$, such that the following
two conditions hold:

\noindent 1) If $\db':d=d'_1+...+d'_{m'}$ is a refinement of $\db$, then
the isomorphism
$$\K^{\db}|_{X^{(\db')}_{disj}}\simeq \K^{\db'}$$
is compatible with the $(\check G)^{\times m}$-actions via the
diagonal map $(\check G_\bF)^{\times m}\to (\check G_\bF)^{\times m'}$.

\noindent 2) If $\db':d=d'_1+...+d'_{m'}$ is obtained from $\db$ by
a permutation, then the isomorphism of perverse sheaves
induced by the isomorphism $X^{\db}\to X^{\db'}$ is compatible with the
$(\check G_\bF)^{\times m}$-actions.

\medskip

Note that for $d=d_1+d_2$, the functor $\K_1,\K_2\mapsto sum_{d_1,d_2}(\K_1\boxtimes \K_2)$
gives rise to a functor
\begin{equation}  \label{conv of rep}
\star:P^{\check G,d_1}\times P^{\check G,d_2}\to P^{\check G,d}
\end{equation}

The next result follows formally from \thmref{satake}:

\begin{thm} \label{extended satake}
For every $d$ there is a canonical equivalence of categories
$P^{\CG^d}(\Gr^d_G)\to P^{\check G,d}$, compatible with
the functors $\star$ of \eqref{conv of sheaves} and \eqref{conv of rep}.
\end{thm}
The proof will be given in Appendix B.

We will denote by $P^{\check G,d}_\bO$ the corresponding
version of $P^{\check G,d}$ with $\bO$-coefficients. \thmref{extended satake}
remains valid in this context as well, i.e.,
$P^{\check G,d}_\bO$ is equivalent to the category
$P^{\CG^d}(\Gr^d_G)_\bO$ of $\CG$-equivariant perverse 
sheaves with $\bO$-coefficients on $\Gr^d_G$.

\ssec{}   \label{Hecke stack}

Let us take $d=1$ and denote the corresponding ind-scheme $\Gr^1_G$ by $\Gr_{G,X}$,
and the corresponding group-scheme $\CG^1$ simply by $\CG$.
Note that $\Gr_{G,X}$ is just the relative over $X$ version of $\Gr_{G,x}$, i.e., we have
a projection $s:\Gr_{G,X}\to X$ and its fiber over $x\in X$ is $\Gr_{G,x}$. 

It follows from \corref{virasoro}, that to an object $V\in \on{Rep}(\check G_\bF)$ one
can canonically attach a perverse sheaf $\T_{V,X}\in P^{\CG}(\Gr_{G,X})$. In terms of
the equivalence of \thmref{extended satake}, $\T_{V,X}$ corresponds to
the constant sheaf $F_X\otimes V[1]$, as an object of $P^{\check G,1}$.

\medskip

Let $\Bun_G$ be the moduli stack of principal $G$-bundles on $X$.
Let us recall the definition of the Hecke functors
$H:\on{Rep}(\check G_\bF)\times D(\Bun_G)\to D(\Bun_G\times X)$.

Let $\H_{G,X}$ be the Hecke stack, i.e, the relative over $\Bun_G$ version of 
$\Gr_{G,X}$. More precisely, $\H_{G,X}$ classifies the data of $(x,\P_G,\P'_G,\beta)$,
where $x\in X$, $\P_G,\P'_G$ are principal $G$-bundles on $X$, and $\beta$ is
an isomorphism $\P_G|_{X-x}\simeq \P'_G|_{X-x}$. We will denote by
$\hl$ (resp., $\hr$) the natural map of stacks $\H_{G,X}\to \Bun_G$ that remembers the
data of $\P'_G$ (resp., $\P_G$). We will view $\H_{G,X}$ as a fibration 
over $\Bun_G$ via $\hl$, with the typical fiber $\Gr_{G,X}$:

$$\Bun_G \overset{\hl}\longleftarrow \H_{G,X} \overset{\hr\times s}\longrightarrow
\Bun_G\times X.$$

Due to the $\CG$-equivariance condition, to every $\S\in D(\Bun_G)$
and $\T\in P^{\CG}(\Gr_{G,X})$ we can associate their twisted external product
$\T\tboxtimes \S\in D(\H_{G,X})$. We define the Hecke functor
$$H(V,\S):=(\hr\times s)_!(\T_{V,X}\tboxtimes \S)\in D(\Bun_G\times X).$$

More generally, if $V_1,...,V_d$ is a collection of objects of $\on{Rep}(\check G_\bF)$,
by iterating the above construction, for $\S\in D(\Bun_G)$ we obtain an object
$$H(V_1\boxtimes...\boxtimes V_d,\S)\in D(\Bun_G\times X^d).$$

As in \cite{Ga}, one shows that for any $\S$ as above, 
$H(V_1\boxtimes...\boxtimes V_d,\S)$ is ULA with respect to the projection
$\Bun_G\times X^d\to X^d$.

\begin{prop}  \label{hecke action}
Let $V_1,V_2$ be two objects of $\on{Rep}(\check G_\bF)$. 

\smallskip

\noindent {\em (1)}
Let $\sigma$ be the transposition acting on $X\times X$. We have a functorial
isomorphism 
$$\sigma^*(H(V_1\boxtimes V_2,\S))\simeq H(V_2\boxtimes V_1,\S),$$
whose square is the identity map.

\smallskip

\noindent {\em (2)}
The restriction $H(V_1\boxtimes V_2,\S)|_{\Bun_G\times \Delta(X)}$
identifies canonically with $H(V_1\otimes V_2,\S)[1]$.

\end{prop}

The above proposition allows us to introduce Hecke eigen-sheaves.
Let $E_{\check G}$ be a $\check G_\bF$-local system on $X$, viewed
as a tensor functor $V\mapsto E^V_{\check G}$ from $\on{Rep}(\check G_\bF)$
to the category of $\bF$-local systems on $X$.

We say that $\S_{E_{\check G}}\in D(\Bun_G)$ is a Hecke eigen-sheaf with respect to
$E_{\check G}$ if we are given an isomorphism of functors
$\on{Rep}(\check G_\bF)\to D(\Bun_G)$:
\begin{equation} \label{eigencondition}
\alpha(V):H(V,\S_{E_{\check G}})\simeq \S_{E_{\check G}}\boxtimes E^V_{\check G}[1],
\end{equation}
such that the conditions 1) and 2) below are satisfied.

Before formulating them, note that for any collection $V_1,...,V_d$ of
objects of $\on{Rep}(\check G_\bF)$, by iterating $\alpha(\cdot)$ we obtain as isomorphism
$$H(V_1\boxtimes...\boxtimes V_d,\S_{E_{\check G}})\simeq 
\S_{E_{\check G}}\boxtimes E^{V_1}_{\check G}[1]\boxtimes...\boxtimes E^{V_m}_{\check G}[1].$$

We require that for $V_1,V_2\in \on{Rep}(\check G_\bF)$,
the following two diagrams commute:

\smallskip

\noindent 1) 
$$
\CD
H(V_1\boxtimes V_2,\S_{E_{\check G}}) @>>> \S_{E_{\check G}}\boxtimes E^{V_1}_{\check G}[1]
\boxtimes E^{V_2}_{\check G}[1] \\
@VVV   @VVV  \\
\sigma^*(H(V_2\boxtimes V_1,\S_{E_{\check G}})) @>>> 
\S_{E_{\check G}}\boxtimes \sigma^*(E^{V_2}_{\check G}[1]\boxtimes 
E^{V_1}_{\check G}[1])
\endCD
$$

\smallskip

\noindent 2) 
$$
\CD
H(V_1\boxtimes V_2,\S)|_{\Bun_G\times \Delta(X)} @>>> \S\boxtimes E^{V_1}_{\check G}[1]\boxtimes 
E^{V_2}_{\check G}[1]|_{\Bun_G\times \Delta(X)} \\
@VVV   @VVV  \\
H(V_1\otimes V_2,\S)[1] @>>> \S\boxtimes E^{V_1\otimes V_2}_{\check G}[1]
\endCD
$$

\ssec{}   \label{generalized Hecke functors}

We need to introduce one more piece of notation related to the
Hecke action.

Let $\H^d_G$ be the relative (over $\Bun_G$) version of $\Gr^d_G$. 
We have the diagram
$$\Bun_G \overset{\hl}\longleftarrow \H^d_G \overset{\hr\times s}\longrightarrow
\Bun_G\times X^{(d)},$$
and we view $\H^d_G$ as a
fibration over $\Bun_G$ via $\hl$ with the typical fiber $\Gr^d_G$.

As in \secref{Hecke stack}, for an object 
$\T\in P^{\CG^d}(\Gr^d_G)$ and $\S\in D(\Bun_G)$ we can form their
twisted external product $\T\tboxtimes \S\in D(\H^d_G)$.
We define
$$H(\T,\S):=(\hr\times s)_!(\T\tboxtimes \S)\in D(\Bun_G\times X^{(d)}).$$
We set $\T\star \S\in D(\Bun_G)$ to be the direct image of $H(\T,\S)$
under $\Bun_G\times X^{(d)}\to \Bun_G$.

By construction, for $\T_1\in P^{\CG^{d_1}}(\Gr^{d_1}_G)$ and
$\T_2\in P^{\CG^{d_2}}(\Gr^{d_2}_G)$, we have a functorial isomorphism
$$(\T_1\star \T_2)\star \S\simeq \T_1\star (\T_2\star \S).$$

We will denote by the same symbol $\star$ the resulting action of 
$P^{\check G,d}$ on $D(\Bun_G)$.

\section{Geometric Langlands conjecture}   \label{intr conj}

In this section will work with $\bF$-sheaves, unless specified otherwise.

\ssec{}   \label{formulation of theorems}

From now on we will specialize to the case $G=GL_n$. In this case
we will denote $\Bun_G$ by $\Bun_n$, and think of principal $G$-bundles
as of rank $n$ vector bundles on $X$. 

Let $V_0$ denote the tautological $n$-dimensional representation of $\check G=GL_n$.
A $\check G$-local system on $X$ is the same thing as a $n$-dimensional local 
system, the correspondence being $E_{\check G}\mapsto E:=E^{V_0}_{E_{\check G}}$.

Here is the formulation of the geometric Langlands conjecture:

\begin{conj} \label{geometric langlands}
If $E$ is an absolutely irreducible $n$-dimensional local system on $X$,
then there exists a perverse sheaf $\S_E\in P(\Bun_n)$, which
is a Hecke eigen-sheaf with respect to $E$. Moreover, $\S_E$
is cuspidal and irreducible on every connected component of $\Bun_n$.
\end{conj}

\medskip

If $\S_E\in D(\Bun_n)$ and $E$ is an $n$-dimensional local system
one can formulate a condition, weaker than the Hecke eigen-property:
We say that $\S_E$ is a weak Hecke eigen-sheaf with respect to $E$,
if we are given isomorphisms $\alpha(\cdot)$ as in \eqref{eigencondition}
for $V$ of the form $\Lambda^i(V_0)$, $i=1,...,n$. 

The weak form of the geometric Langlands conjecture says that if
$E$ is an irreducible local system on $X$, then there exists a 
perverse sheaf $\S_E\in P(\Bun_n)$, which is a {\it weak} Hecke 
eigen-sheaf with respect to $E$, such that $\S_E$ is cuspidal 
and irreducible on every connected component of $\Bun_n$. 

Evidently, the existence of a weak Hecke eigen-sheaf corresponding to
$E$ is sufficient to guarantee the existence of an $\bF$-valued 
cuspidal Hecke eigen-form corresponding to (the $\pi_1$-representation,
corresponding to $E$) in the sense of \conjref{Langlands}.

\medskip

Yet another form of the Hecke eigen-condition, this time specific to
$GL_n$ is as follows:

\smallskip

For an $n$-dimensional local system $E$, we say that $\S_E\in D(\Bun_n)$
has a $GL_n$-Hecke eigen-property with respect to $E$ if we are given
an isomorphism $\alpha(\cdot)$ {\it only for $V=V_0$}, which satisfies
condition 1) of the definition of Hecke eigen-sheaves for
$V_1\simeq V_2\simeq V$.

\begin{conj}   \label{GL_n Hecke}
Let $E$ be an arbitrary $n$-dimensional local system. Then if
$\S_E\in D(\Bun_n)$ has a $GL_n$-Hecke eigen-property with respect to $E$,
then it satisfies, in fact, the full Hecke property.
\end{conj}

When $\on{char}(\bF)=0$, the above conjecture was essentially proved in \cite{Ga2},
using Springer correspondence. We do not have any real evidence in favor
of this conjecture when $\on{char}(\bF)\neq 0$ (except in the case when
$\S_E$ is cuspidal and perverse, which we consider in Appendix A).
However, we have the following assertion:

\begin{lem}  \label{other Hecke}
Assume that $\on{char}(\bF)>n$, and let $\S_E\in D(\Bun_n)$ have a
$GL_n$-Hecke  eigen-property with respect to $E$. Then $\S_E$ is a weak 
Hecke eigen-sheaf with respect to $E$. 
\end{lem}

The proof of this lemma given in \cite{FGV} in the case when
$\on{char}(\bF)=0$ is applicable here, since the proof relies on the
semi-simplicity of representations of $\Sigma_i$ on $\bF$-vector
spaces, which is valid if $i<\on{char}(\bF)$.

\medskip

We will prove the following:

\begin{thm} \label{main}
Assume that $\on{char}(\bF)>2n$. Then for every irreducible local system $E$ 
on $X$ there exists a (cuspidal, irreducible on every connected component)
perverse sheaf $\S_E\in P(\Bun_n)$, which
is a $GL_n$-Hecke eigen-sheaf with respect to $E$. 
\end{thm}

Combined with \lemref{other Hecke} this proves the weak
form of the geometric Langlands conjecture, and hence 
\conjref{Langlands}, assuming that $\on{char}(\bF)>2n$.

\medskip

In fact, by relying on some more unpublished work, one can strengthen this
result, and prove the following:

\begin{thm}   \label{strong form}
Assume that $\on{char}(\bF)\neq 2$. Then for an absolutely irreducible local system $E$ 
on $X$ there exists a (cuspidal, irreducible on every connected component)
perverse sheaf $\S_E\in P(\Bun_n)$, which
is a Hecke eigen-sheaf with respect to $E$. 
\end{thm}

We will sketch the proof of this theorem in Appendix A.
Of course, \thmref{strong form} implies de Jong's conjecture for any $l\neq 2$.

\ssec{}

To prove \thmref{main} we will follow the strategy of \cite{FGV}.
Let $\Mod^d_n$ be the stack of upper modifications of length $d$;
note that $\Mod^d_n$ is a closed substack of $\H^d_{GL_n}$ corresponding
to the condition that the generic isomorphism of bundles $\beta:\M\to \M'$
is such that it extends to a regular map of coherent sheaves.

We set $\T^d_E$ to be the perverse sheaf on $\H^d_{GL_n}$ corresponding 
via the equivalence of \thmref{extended satake} to the following object of
$P^{d,GL_n}$:
$$(E\otimes V_0)^{(d)}[d]$$
(cf. \secref{symmetric powers} below, where symmetric powers of local systems
are discussed). One easily shows that $\T^d_E$ is indeed supported on $\Mod^d_n$.

Let $\on{Coh}^d_0$ be the stack classifying torsion sheaves of length $d$.
There is a natural smooth projection $\Mod^d_n\to \on{Coh}^d_0$, and
$\T^d_E$ is isomorphic (up to a shift) to the pull-back of a canonical perverse
sheaf, denoted $\L^d_E$ on $\on{Coh}^d_0$, called the Laumon sheaf.

More explicitly, $\L^d_E$
is the Goresky-MacPherson extension of its own restriction to the open
substack $\overset{\circ}{\on{Coh}}{}^d_0$, corresponding to regular semi-simple
torsion sheaves. This restriction is the pull-back of
$\overset{\circ}E{}^{(d)}$ under the natural smooth morphism
$\overset{\circ}{\on{Coh}}{}^d_0\to \overset{\circ}X{}^{(d)}$.

\medskip

We will need the following statement, whose $\on{char}=0$ version was proved in 
\cite{La}. For $d=d_1+d_2$, let $\on{Fl}^{d_1,d_2}_0$ denote the stack classifying short
exact sequences $0\to J_1\to J\to J_2\to 0$, where $J_1$ and $J_2$
are torsion coherent sheaves of lengths $d_1$ and $d_2$, respectively.

Let $\p$ denote the natural projection $\on{Fl}^{d_1,d_2}_0\to \on{Coh}^d_0$,
and let $\q$ be the projection $\on{Fl}^{d_1,d_2}_0\to \on{Coh}^{d_1}_0\times 
\on{Coh}^{d_2}_0$.

\begin{thm}    \label{Hecke-Laumon property}
For any local system $E$ we have:
$$\q_!\circ \p^*(\L^d_E)\simeq \L^{d_1}_E\boxtimes \L^{d_2}_E.$$
\end{thm}

\ssec{Proof of \thmref{Hecke-Laumon property}}

Consider first the open substack of $\on{Coh}^{d_1}_0\times \on{Coh}^{d_2}_0$
equal to $\q(\p^{-1}(\overset{\circ}{\on{Coh}}{}^d_0))$. Over it both maps
$\q$ and $\p$ are isomorphisms and the isomorphism stated in the theorem is
evident. Therefore, we have to show that $\q_!\circ p^*(\L^d_E)$ is a
perverse sheaf, extended minimally from the above open substack.

The question being local, we can assume that $E$ is trivial. Let us decompose
$E$ as a sum of $1$-dimensional local systems $E=E_1\oplus...\oplus E_n$.
For a partition $\db:d=d^1+...+d^n$, let $\on{Fl}^{\db}_0$ be the 
stack classifying flags of coherent torsion sheaves with successive quotients
having lengths given by $\db$, and let $\p^{\db}$ be the natural projection
$\on{Fl}^{\db}_0\to \on{Coh}^d_0$. Let $\q^{\db}$ denote the map
$\on{Fl}^{\db}_0\to \underset{i}\Pi\, X^{(d^i)}$ obtained by taking supports
of the successive quotients.

From \thmref{extended satake} and \lemref{symmetric power direct sum}, it follows that
$$\L^d_E\simeq \underset{\db}\oplus\, \p^{\db}_!\circ 
(\q^{\db})^*(E_1^{(d^1)}\boxtimes...\boxtimes E_n^{(d^n)}).$$
Note that each $E_i^{(d^i)}$ is the constant sheaf on $X^{(d^i)}$, since $E_1$ 
is one dimensional and trivial.

Thus, we have to compute
$$\q_!\circ \p^*\circ \p^{\db}_!\circ (\q^{\db})^*(\bF_{\on{Fl}^{\db}_0}).$$
Let $Z$ denote the fiber product
$\on{Fl}^{\db}_0\underset{\on{Coh}^d_0}\times \on{Fl}^{d_1,d_2}_0$. It can
be naturally subdivided into locally closed subvarieties (Schubert cells)
numbered by the set 
$$(\Sigma_{d_1}\times \Sigma_{d_2})
\backslash\Sigma_d/(\Sigma_{d^1}\times...\times \Sigma_{d^n}),$$
or, which is the same, of the ways to partition $d_1=d_1^1+...+d_1^n$,
$d_2=d_2^1+...+d_2^n$ with $d_1^i+d_2^i=d^i$. For each such pair of partitions,
let us denote by $Z^{\db_1,\db_2}$ the corresponding subvariety in $Z$.

It is sufficient to show that the direct image of the constant sheaf on $Z^{\db_1,\db_2}$
under
$$Z^{\db_1,\db_2}\hookrightarrow Z\to \on{Fl}^{d_1,d_2}_0 \overset{\q}\to
\on{Coh}^{d_1}_0\times \on{Coh}^{d_2}_0$$ is a perverse sheaf, minimally extended
from the open substack $\q(\p^{-1}(\overset{\circ}{\on{Coh}}{}^d_0))$.

The above map factors as
$$Z^{\db_1,\db_2}\to \on{Fl}^{\db_1}_0\times \on{Fl}^{\db_1}_0
\overset{\p^{\db_1}\times \p^{\db_2}}
\longrightarrow \on{Coh}^{d_1}_0\times \on{Coh}^{d_2}_0,$$
where the first arrow is a generalized smooth fibration into affine spaces,
of relative dimension $0$. Therefore, the resulting object of $D(\on{Coh}^{d_1}_0\times
\on{Coh}^{d_2}_0)$ is isomorphic to 
$\p^{\db_1}_!(\bF_{\on{Fl}^{\db_1}_0})\boxtimes \p^{\db_2}_!(\bF_{\on{Fl}^{\db_2}_0})$,
and the latter is known (cf. \thmref{extended satake}) to be a perverse sheaf, 
minimally extended from the required open locus.

\ssec{}

Proceeding as in \cite{FGV}, we produce from $\T^d_E$ the Whittaker sheaf,
and subsequently an object $\S'_E\in D(\Bun'_n)$, where $\Bun'_n$ is the
stack classifying pairs $(\M\in \Bun_n,\kappa:\Omega^{n-1}\to \M)$.

Let $U\subset \Bun_n$ be the open substack, defined by the condition that
$\M\in U$ if $\on{Ext}^1(\Omega^{n-1},\M)=0$. Let $U'$ be the preimage of $U$
in $\Bun'_n$. Clearly, the restriction of the projection $\pi:\Bun'_n\to \Bun_n$
to $U$ is smooth.

We claim that the restriction of $\S'_E$ to $U'$ is perverse and irreducible on
every connected component. This follows as in \cite{FGV}, using 
\thmref{Hecke-Laumon property} from the vanishing result, \thmref{vanishing},
discussed below.

Having established \thmref{vanishing}, and hence perversity and irreducibility
of $\S'_E$, the next step is to show that $\S'_E$ descends to a perverse sheaf
on $\Bun_n$. We do it by the same argument involving Euler-Poincar\'e 
characteristics as in {\it loc.cit.} Namely, we have to show that the
Euler-Poincar\'e characteristics of $\S'_E$ are constant along the fibers of
the projection $\pi$.

As in \cite{FGV}, using Deligne's theorem, we show that the 
Euler-Poincar\'e characteristics of $\S'_E$ are independent of the local system.
Therefore, it suffices to show that the Euler-Poincar\'e characteristics of 
$\S'_{E_0}$ are constant along the fibers of $\pi$, where $E_0$ is the {\it trivial}
$n$-dimensional local system. When we work with $\bF'$-sheaves
(where $\bF'$ is a local field of characteristic 0 with residue field $\BF_l$),
rather than with
$\bF=\BF_l((t))$-sheaves, the corresponding fact follows from Sect. 6 of \cite{FGV}. 
Therefore, it suffices to show that the Euler-Poincar\'e characteristics of $\S'_{E_0}$ at
a given point of $\Bun'_n$ are the same in the $\bF'$ and $\bF$-situations. We will
show this by comparing both sides with $\BF_l$-sheaves.

\medskip

Thus, let $\bF_0$ be any of the local fields ($\bF$ or $\bF'$), and let
$\bO_0$ be the corresponding local ring. By a subscript ($\bF_0$, $\bO_0$ or $\BF_l$)
we will indicate which of the sheaf-theoretic contexts we are working in.

Consider the corresponding category $P^{d,GL_n}_{\bO_0}$. We can form the object
$(E_0\otimes V_0)^{(d)}_{\bO_0}\in P^{d,GL_n}_{\bO_0}$, which under $\bO_0\to \bF_0$ and
$\bO_0\to \BF_l$ specializes to $(E_0\otimes V_0)^{(d)}_{\bF_0}$ and
$(E_0\otimes V_0)^{(d)}_{\BF_l}$, respectively. From \secref{symmetric powers}
it follows that $(E_0\otimes V_0)^{(d)}_{\bO_0}$ is $\bO_0$-flat.

Using \thmref{extended satake} for $\bO_0$,
from $(E_0\otimes V_0)^{(d)}_{\bO_0}$ we produce the corresponding
$\bO_0$-perverse sheaf on $\Mod^d_n$, which is also $\bO_0$-flat. Finally, we produce
the complex of $\bO_0$-sheaves $(\S'_{E_0})_{\bO_0}$ on $\Bun'_n$. By the above flatness
property
$$(\S'_{E_0})_{\bO_0}\overset{L}{\underset{\bO_0}\otimes} \BF_l\simeq
(\S'_{E_0})_{\BF_l} \text{ and }
(\S'_{E_0})_{\bO_0}\overset{L}{\underset{\bO_0}\otimes} \bF_0\simeq
(\S'_{E_0})_{\bF_0},$$ 
where $(\S'_{E_0})_{\BF_l}$ and $(\S'_{E_0})_{\bF_0}$ are the corresponding
complexes of $\BF_l$-sheaves and $\bF_0$-sheaves, respectively, on $\Bun'_n$.

Let $\K_{\bO_0}$ (resp., $\K_{\BF_l}$, $\K_{\bF_0}$) be the fiber
of $(\S'_{E_0})_{\bO_0}$ (resp., $(\S'_{E_0})_{\BF_l}$, $(\S'_{E_0})_{\bF_0}$)
at a given point of $\Bun'_n$.
This is an object of the bounded derived category of finitely-generated
$\bO_0$-modules (resp., finite-dimensional $\BF_l$- or $\bF_0$-vector spaces), 
and we still have
$$\K_{\bO_0}\overset{L}{\underset{\bO_0}\otimes} \BF_l\simeq
\K_{\BF_l} \text{ and }
\K_{\bO_0}\overset{L}{\underset{\bO_0}\otimes} \bF_0\simeq
\K_{\bF_0}.$$
Under these circumstances we always have: $\chi(\K_{\bF_0})=\chi(\K_{\BF_l})$.
Indeed, it suffices to consider separately the cases when $\K_{\bO_0}$ is a flat
finitely-generated $\bO_0$-module (in which case the assertion is evident),
or when $\K_{\bO_0}$ is torsion. In the latter case, we can assume that 
$\K_{\bO_0}\simeq \BF_l$. Then $\K_{\bF_0}=0$, and $\K_{\BF_l}$
has $1$-dimensional cohomologies in degrees $-1$ and $0$, i.e., 
$\chi(\K_{\BF_l})=0$.

\medskip

This proves the fact that $\S'_E$ descends to a perverse sheaf 
$\S_E$ on $\Bun_n$. The fact that $\S_E$ satisfies the Hecke property
follows from \thmref{Hecke-Laumon property} as in \cite{FGV}, Section 8.
The cuspidality of $\S_E$ also follows from the vanishing result, \thmref{vanishing}.

\section{The vanishing result}   \label{sect vanishing}

\ssec{} \label{interpretation of Laumon}

To state \thmref{vanishing} we have to recall the definition of
the averaging functor
$\on{Av}^d_E:D(\Bun_n)\to D(\Bun_n)$. It is defined for $d\in {\mathbb N}$ 
and a local system $E$ of an arbitrary rank.

By definition,
$$\on{Av}^d_E(\S)=\hr_!(\hl{}^*(\S)\otimes \T^d_E).$$
In other words,
$$\on{Av}^d_E(\S)=(E\otimes V_0)^{(d)}[d]\star \S,$$
in the notation of \secref{generalized Hecke functors}.

The key step in the proof of perversity and irreducibility of
$\S'_E$ on $U'$ is the following:

\begin{thm}  \label{vanishing}
Let $E$ be absolutely irreducible, of rank $m$ with $m>n$ 
and $d$ be $>(2g-2)\cdot m\cdot n$. Then the functor
$\on{Av}^d_E$ is identically equal to $0$.
\end{thm}

\ssec{}

To prove \thmref{vanishing} we will have to analyze separately the
cases of $\on{char}(k)=0$ and $k$ of positive characteristic.
Let us show that the former case reduces to the latter.
(Of course, for de Jong's conjecture we need 
the case $k=\ol{\BF}_q$.)

Indeed, if $k$ is of characteristic $0$, we can replace it
by ${\mathbb C}$ and work with sheaves in the analytic topology. 
In this case, we can consider sheaves with coefficients in
an arbitrary field $F$ of characteristic $l$. If $F$ is 
finite, the local system $E$ is defined over a finitely
generated sub-field of $k$, and the standard procedure
reduces us to the case of a finite ground. 

The case of a general $F$ reduces to that of $\BF_l$. 
Indeed, since the fundamental group of a curve is finitely
generated, we can assume that our local system has 
coefficients in a ring $A$, finitely generated over $\BF_l$. 
Then the vanishing of the functor $\on{Av}^d_E$ over $A$ 
would follow from the corresponding assertion at all geometric
points of $A$, whose residue fields are finite.

\ssec{}

Thus, for the rest of this section we will assume that
the ground field $k$ is of positive characteristic of
order prime to $l$. We need this assumption in order to
have the Artin-Schreier sheaf on the affine line, and
the Fourier transform, which is used in the 
definition of quotient categories $\wt{D}(\Bun_n)$, 
see below.

The proof of \thmref{vanishing} will 
follow the same lines as the proof of
the analogous statement in the situation of $\on{char}=0$ coefficients
in \cite{Ga}. We have:

\begin{thm}  \label{exactness}
Under the assumptions of \thmref{vanishing},
the functor $\on{Av}^d_E$ is exact in the sense of the perverse
t-structure on $D(\Bun_n)$.
\end{thm}

We show as in \cite{Ga}, Appendix (or, alternatively, as in Sect. 2.1 of 
{\it loc. cit.}, which is slightly more cumbersome) that
\thmref{exactness} implies \thmref{vanishing}. Thus, our goal from now
on is to prove \thmref{exactness}.

The first step is to introduce a quotient triangulated category
$\wt{D}(\Bun_n)$ of $D(\Bun_n)$, which has Properties 0,1,2 of
\cite{Ga}, Sect. 2.12. The construction of $\wt{D}(\Bun_n)$ and
the verification of its properties given in Sect.4-8 of {\it loc.cit.}
goes through without modification in our situation.

The next step is to prove an analog of Theorem 2.14 of {\it loc. cit.}, that says that
the functor $\on{Av}^1_E$ is exact on the quotient category $\wt{D}(\Bun_n)$,
provided that $E$ is absolutely irreducible of rank strictly greater than $n$.
Again, the argument presented in {\it loc. cit.} is applicable, since
it only involves the action of symmetric groups $\Sigma_k$ with $k\leq 2n$,
and we have made the assumption $\on{char}(\bF)>2n$.

\medskip

The final step, which will require some substantial modifications in the
case of coefficients of positive characteristic is the following:

\begin{prop}   \label{functor on quotient}
Let $\wt{D}(\Bun_n)$ be a triangulated quotient category of $D(\Bun_n)$,
satisfying Properties 0 and 1 above. 

\smallskip

\noindent{\em (1)}
For any $E$, the functor $\on{Av}^d_E$ descends to a
well-defined functor on $\wt{D}(\Bun_n)$.

\smallskip

\noindent{\em (2)}
If $E$ is such that the functor $\on{Av}^1_E$ is exact on $\wt{D}(\Bun_n)$,
then so is $\on{Av}^d_E$ for any $d$.
\end{prop}

Once \propref{functor on quotient} is proved, we finish the proof
of \thmref{exactness} as in Sect. 2.16 of \cite{Ga}.

\ssec{}    \label{action on quotient}

Let $\wt{D}(\Bun_n\times S)$ be the system of quotient categories 
of $D(\Bun_n\times S)$, satisfying Properties 0 and 1. Recall that
the generalized Hecke functors $H(\cdot,\cdot)$ of 
\secref{generalized Hecke functors}. We will prove the following:

\begin{prop}  \label{Hecke action on quotient}
For any $\T\in P^{\CG^d}(\Gr^d_G)$ the functor
$$D(\Bun_n)\to D(\Bun_n\times X^{(d)}): \S\mapsto H(\T,\S)$$
descends to a well-defined exact functor 
$\wt{D}(\Bun_n)\to \wt{D}(\Bun_n\times X^{(d)})$.
\end{prop}

Clearly, \propref{Hecke action on quotient} implies 
\propref{functor on quotient}(1). In addition, we have the following
corollary:

\begin{cor}  \label{dimension}
Let $\T$ be an object of $P^{\CG^d}(\Gr^d_G)$ supported 
over a subvariety of dimension $\leq i$ of $X^{(d)}$.
Then the functor $\S\mapsto \T\star \S$ has the 
cohomological amplitude at most $[-i,i]$ on $\wt{D}(\Bun_n)$.
\end{cor}

The rest of this subsection is devoted to the proof of
\propref{Hecke action on quotient}.

\medskip

\begin{lem}   \label{aux hecke}
The functor 
$$H(V,\cdot):D(\Bun_n)\to D(\Bun_n\times X)$$
(cf. \secref{Hecke stack}) descends to a well-defined exact functor 
$\wt{D}(\Bun_n)\to \wt{D}(\Bun_n\times X)$ for any $V\in \on{Rep}(\check G_\bF)$.
\end{lem}

\begin{proof}

First, we claim that this is true for $V=\Lambda^i(V_0)$. This follows
from the Springer correspondence (applicable here, since $\on{char}(\bF)>n$)
as in Prop. 1.11 of \cite{Ga}. The fact that $H(V,\cdot)$ is well-defined
for any $V$ follows now, since the classes of the
representations of the form $\Lambda^i(V_0)$ generate the Grothendieck
ring of $\on{Rep}(\check G_\bF)$.

To prove the exactness assertion, it is enough to assume that $V$ is 
irreducible. We will proceed by induction on the length of the highest 
weight of $V$. Since the statement is essentially Verdier self-dual,
it is enough to prove that the functor 
$H(V,\cdot):\wt{D}(\Bun_n)\to \wt{D}(\Bun_n\times X)$ is right-exact.
However, for any such $V$, there exists a representation $V'$ isomorphic
to a tensor product of representations of the form $\Lambda^i(V_0)$,
together with a surjection $V'\twoheadrightarrow V$, such that its
kernel, $V''$, is an extension of irreducible representations with smaller
highest weights. 

For $\S\in \wt{P}(\Bun_n)$ we have a long exact cohomology sequence
$$...\to h^i(H(V'',\S))\to h^i(H(V',\S))\to h^i(H(V,\S))\to h^{i+1}(H(V'',\S))\to...$$
and by the induction hypothesis, we conclude that $H(V,\cdot)$ is right exact.

\end{proof}

Consider the diagonal stratification of $X^{(d)}$ numbered by the partitions
$\db:d=d_1+...+d_k$. For each such partition consider the space
$$X_{\db}:=(\underset{k \text{ times }}{\underbrace{X\times ...\times X}})_{disj},$$
which covers in a finite and \'etale way the corresponding stratum in $X^{(d)}$.

Let us denote by $\Gr_{G,\db}$ the fiber product
$X_{\db}\underset{X^{(d)}}\times \Gr^d_G$. Note that $\Gr_{G,\db}$ is isomorphic
to $$\underset{k \text{ times }}{\underbrace{\Gr_{G,X}\times...\times \Gr_{G,X}}}\underset{X^k}
\times X_{\db}.$$
We can consider the group scheme $\CG_{\db}:=\CG^{\times k}|_{X_{\db}}$, the
category $P^{\CG_{\db}}(\Gr_{G,\db})$
of $\CG_{\db}$-equivariant perverse sheaves on $\Gr_{G,\db}$, and the corresponding
Hecke functor
$$H(\cdot,\cdot):P^{\CG_{\db}}(\Gr_{G,\db})\times 
D(\Bun_n)\to D(\Bun_n\times X_{\db}).$$

To prove the proposition, it is enough to show that the latter functor
descends to a well-defined exact functor
$P^{\CG_{\db}}(\Gr_{G,\db})\times \wt{D}(\Bun_n)\to \wt{D}(\Bun_n\times X_{\db})$.

\medskip

Note that every irreducible object of $P^{\CG_{\db}}(\Gr_{G,\db})$ has the form
$$(\T_{V_1,X}\boxtimes...\boxtimes \T_{V_k,X}|_{X_{\db}})\otimes \K,$$
where $\K$ is an irreducible perverse sheaf on $X_{\db}$, and 
$V_1,...,V_k$ are irreducible objects of $\on{Rep}(\check G_\bF)$ 
(the notation $\T_{V,X}$ is as in \secref{Hecke stack}).

For such an object of $P^{\CG^d}(\Gr^d_G)$ the above functor $H(\cdot,\cdot)$
takes the form
$$\S\mapsto H(V_1\boxtimes...\boxtimes V_k,\S)|_{X_{\db}}\otimes \K.$$
From \lemref{aux hecke} we obtain that this functor indeed descends to
$\wt{D}(\Bun_n)$. The exactness statement also follows from \lemref{aux hecke},
since any object of the form $$H(V_1\boxtimes...\boxtimes V_k,\S)\in D(\Bun_n\times X^k)$$
is ULA with respect to the projection $\Bun_n\times X^k\to X^k$, cf. Lemma 3.7 of 
\cite{Ga}.

\ssec{}

Now we are ready to prove \propref{functor on quotient}(2). For that we need
one more piece of preparatory material, namely, the notion of external
exterior power of a local system. This notion is discussed in \secref{sect ext ext}.
Thus, for any integer $d$ we have a perverse sheaf
$\Lambda^{!,(d)}(E\otimes V_0)$ on $X^{(d)}$. Moreover, by 
\thmref{external exterior}(2), $\Lambda^{!,(d)}(E\otimes V_0)$ is
naturally an object of $P^{d,GL_n}$.

\medskip

According to \propref{Hecke action on quotient}, we have a well-defined functor
$$\star:P^{d,GL_n}\times \wt{D}(\Bun_n)\to \wt{D}(\Bun_n),$$
and we must prove the exactness of $(E\otimes V_0)^{(d)}\star ?$.

Assume that $\S$ belongs to $\wt{P}(\Bun_n)$. We will prove by induction that 
$(E\otimes V_0)^{(d)}\star \S$ also belongs to $\wt{P}(\Bun_n)$. Thus, we assume
that the statement holds for $d'<d$. Since the situation is essentially Verdier
self-dual, it is enough to show that $(E\otimes V_0)^{(d)}\star \S\in \wt{D}^{\leq 0}(\Bun_n)$.

Consider the complex 
\begin{align*}
&\Lambda^{!,(d)}(E\otimes V_0)\to \Lambda^{!,(d-1)}(E\otimes V_0)\star (E\otimes V_0)\to...\to
\Lambda^{!,(i)}(E\otimes V_0)\star (E\otimes V_0)^{(d-i)}\to \\
&\Lambda^{!,(i-1)}(E\otimes V_0)\star (E\otimes V_0)^{(d-i+1)}\to...\to
(E\otimes V_0)\star (E\otimes V_0)^{(d-1)}\to (E\otimes V_0)^{(d)}
\end{align*}
of objects of $P^{d,GL_n}$, given by \thmref{external exterior}. Since this
complex is exact, it is enough to show that each
$$(\Lambda^{!,(i)}(E\otimes V_0)\star (E\otimes V_0)^{(d-i)})\star \S$$
belongs to $\wt{D}^{\leq i-1}(\Bun_n)$ for $i=1,...,d$.

By the induction hypothesis, we know that 
$(E\otimes V_0)^{(d-i)}\star \S\in \wt{D}^{\leq 0}(\Bun_n)$. Therefore, it suffices
to show that the functor
$$\S\mapsto \Lambda^{!,(i)}(E\otimes V_0)\star \S$$
sends $\wt{D}^{\leq 0}(\Bun_n)$ to $\wt{D}^{\leq i-1}(\Bun_n)$

Note that by \corref{dimension}, for $\S\in \wt{D}^{\leq 0}(\Bun_n)$, 
the object $\Lambda^{!,(i)}(E\otimes V_0)\star \S$ does belong to 
$\wt{D}^{\leq i}(\Bun_n)$, since $X^{(i)}$ is $i$-dimensional. Therefore, it
suffices to show that the top cohomology
$h^i(\Lambda^{!,(i)}(E\otimes V_0)\star \S)$ vanishes.

\medskip

Consider the object $j_{!*}(\overset{\circ}\Lambda{}^{!,(i)}(E\otimes V_0))[i]\in P^{i,GL_n}$.
We have an injective map 
$$j_{!*}(\overset{\circ}\Lambda{}^{!,(i)}(E\otimes V_0)[i])\hookrightarrow
\Lambda^{!,(i)}(E\otimes V_0),$$
and the cokernel is an object of $P^{i,GL_n}$ supported on a subvariety of dimension
$\leq i$. Therefore, by \corref{dimension} and the long exact sequence,
it suffices to show that
$$h^0\left(j_{!*}(\overset{\circ}\Lambda{}^{!,(i)}(E\otimes V_0)[i])\star  \S\right)=0$$
for $\S\in \wt{D}^{\leq 0}(\Bun_n)$.
We have a surjection of perverse sheaves on $\overset{\circ}X{}^{(i)}$
$$j^*((\sym_i)_!((E\otimes V_0)^{\boxtimes i}))[i]\to 
\overset{\circ}\Lambda{}^{!,(i)}(E\otimes V_0)[i],$$ and hence also a surjection
$$(\sym_i)_!((E\otimes V_0)^{\boxtimes i})[i]\twoheadrightarrow 
j_{!*}(\overset{\circ}\Lambda{}^{!,(i)}(E\otimes V_0)[i])$$
of objects of $P^{i,GL_n}$.
Again, by \corref{dimension} and the long exact sequence, it suffices to show
that $$h^0\left((\sym_i)_!((E\otimes V_0)^{\boxtimes i})[i]\star  \S\right)=0$$
for $\S\in \wt{D}^{\leq 0}(\Bun_n)$. However,
$$(\sym_i)_!((E\otimes V_0)^{\boxtimes i})[i]\star  \S=
\underset{i \text{ times}}{\underbrace{\wt{\on{Av}}^1_E\circ ...\circ \wt{\on{Av}}^1_E}}(\S)[i].$$
By the assumption, the functor $\wt{\on{Av}}^1_E$ is exact. Hence, the above
expression has zero cohomologies in all the degrees $>-i$, and in particular, in degree $0$,
if $\S\in \wt{D}^{\leq 0}(\Bun_n)$.

\section{Symmetric and exterior powers of local systems}  \label{sect powers}

In this section we will work with $\bF$-vector spaces, and $\bF$-local systems on $X$.
However, the same results extend to flat $\bO$-modules and local systems.

\ssec{}   \label{usual exterior powers}

Let $V$ be a vector space over $\bF$. Let us denote by $\Sym^{!,2}(V)$ the subspace
of $V^{\otimes 2}$ consisting of flip-invariant vectors. Note that
$\Sym^{!,2}(V)$ is spanned by vectors of the form $v\otimes v$, $v\in V$.
Let $\Sym^{!,d}(V)$ be the subspace of $V^{\otimes d}$ 
consisting of invariants of the symmetric group $\Sigma_d$. Of course, for
$d=1$, $\Sym^{!,1}(V)=V$ and for $d\geq 2$,
$$\Sym^{!,d}(V)=\underset{1\leq i\leq d-1}\cap\, V^{\otimes i-1}\otimes \Sym^{!,2}(V)\otimes
V^{\otimes d-1-i},$$
since $\Sigma_d$ is generated by the simple reflections.

Let $\Lambda^{!,2}(V)$ be the subspace of $V^{\otimes 2}$ spanned by
vectors of the form $v\otimes w-w\otimes v$, $v,w\in V$. For $d\geq 2$, let 
$\Lambda^{!,d}(V)$ be the subspace of $V^{\otimes d}$ equal to the intersection
$$\underset{1\leq i\leq d-1}\cap\, V^{\otimes i-1}\otimes \Lambda^{!,2}(V)\otimes
V^{\otimes d-1-i}.$$

Note that when $\on{char}(\bF)\neq 2$, we can define $\Lambda^{!,2}(V)$as 
the subspace of $\sigma$-anti-invariants in $V^{\otimes 2}$, where 
$\sigma$ is the transposition; hence
in this case $\Lambda^{!,d}(V)$ coincides with the subspace of 
$\Sigma_d$-anti-invariants in $V^{\otimes d}$.

Let $\Sym^{*,d}(V)$ be the quotient space of $V^{\otimes d}$ by the sum of subspaces
of the form $V^{\otimes i-1}\otimes \Lambda^{!,2}(V)\otimes V^{\otimes d-1-i}$ for 
$1\leq i\leq d-1$. In other words, $\Sym^{*,d}(V)$ is the space of 
$\Sigma_d$-coinvariants in $V^{\otimes d}$.

Let $\Lambda^{*,d}(V)$ be the quotient of $V^{\otimes d}$ by the 
sum of subspaces of the form 
$V^{\otimes i-1}\otimes \Sym^{!,2}(V)\otimes V^{\otimes d-1-i}$
for $1\leq i\leq d-1$. 

Note that as representations of $GL(V)_\bF$, $\Lambda^{*,d}(V)$ and $\Lambda^{!,d}(V)$
are irreducible and isomorphic to one-another (but, of course, the isomorphism 
is not given by the map $\Lambda^{!,d}(V)\to V^{\otimes d}\to \Lambda^{*,d}(V)$,
as the latter is zero if $\on{char}(\bF)$ divides $d$). The isomorphism in question is induced by the
endomorphism of $V^{\otimes d}$ given by $\underset{\sigma\in \Sigma_d}{\bf \Sigma}\,
sign(\sigma)\cdot \sigma$. We will sometimes use the notation $\Lambda^i(V)$ to
denote either of the above vector spaces.

By definition, $\Lambda^1(V)=V$
and $\Lambda^0(V)=\Sym^{*,0}(V)=\Sym^{!,0}(V)=\bF$.
If $d>\dim(V)$ one easily shows that $\Lambda^d(V)=0$.

\begin{lem}  \label{linear algebra powers}
We have:

\smallskip

\noindent{\em (1)}
If $V$ is finite-dimensional, there are canonical isomorphisms 
$$\Sym^{*,d}(V^{*})\simeq (\Sym^{!,d}(V))^{*} \text{ and }
\Lambda^{*,d}(V^{*})\simeq (\Lambda^{!,d}(V))^{*},$$
where ${}^{*}$ denotes the dual vector space.

\smallskip

\noindent{\em (2)}
If $V\simeq V_1\oplus V_2$, there are canonical isomorphisms
\begin{align*}
&\Sym^{*,d}(V)\simeq \underset{d_1+d_2=d}\oplus\, 
\Sym^{*,d_1}(V_1)\otimes \Sym^{*,d_2}(V_2) \text{ and } \\
&\Lambda^{*,d}(V)\simeq \underset{d_1+d_2=d}\oplus\, 
\Lambda^{*,d_1}(V_1)\otimes \Lambda^{*,d_2}(V_2),
\end{align*}
and similarly for the !-versions.

\smallskip

\noindent{\em (3)}
The natural maps $\Lambda^{!,d}(V)\to \Lambda^{!,d-1}(V)\otimes V$ and
$\Sym^{*,d-1}(V)\otimes V\to \Sym^{*,d}(V)$ give rise to the long exact
sequence (called the Koszul complex)
\begin{align*}
&\Lambda^{!,d}(V)\to \Lambda^{!,d-1}(V)\otimes V\to
...\to \Lambda^{!,i}(V)\otimes \Sym^{*,d-i}(V)\to \\
&\Lambda^{!,i-1}(V)\otimes \Sym^{*,d-i+1}(V)\to...\to  
V\otimes \Sym^{*,d-1}(V)\to \Sym^{*,d}(V),
\end{align*}
and the maps $\Sym^{!,d}(V)\to \Sym^{!,d-1}(V)\otimes V$ and
$\Lambda^{*,d-1}(V)\otimes V\to \Lambda^{*,d}(V)$ give rise to the long exact
sequence 
\begin{align*}
&\Sym^{!,d}(V)\to \Sym^{!,d-1}(V)\otimes V\to
...\to \Sym^{!,i}(V)\otimes \Lambda^{*,d-i}(V)\to \\
&\Sym^{!,i-1}(V)\otimes \Lambda^{*,d-i+1}(V)\to...\to 
V\otimes \Lambda^{*,d-1}(V)\to \Lambda^{*,d}(V).
\end{align*}
Moreover, the above long exact sequences transform to one-another
under duality, if $V$ is finite-dimensional.

\end{lem}

Points (1) and (2) of this lemma are straightforward. For the proof of (3), 
we reduce the assertion to the case of $\dim(V)=1$ using (2).

\ssec{}  \label{symmetric powers}

Let now $X$ be a smooth curve, let $X^{(d)}$ denote its $d$-th symmetric power,
i.e., $X^{(d)}=X^d/\Sigma_d$, and let $\sym_d$ denote the projection
$X^d\to X^{(d)}$. It is well-known that $X^{(d)}$ is smooth, and in particular,
the map $\sym_d$ is flat.

Let $E$ be a local system on $X$; we are going to recall the construction of
its $d$-th external symmetric power. Recall that $sum_{d_1,d_2}$ denotes the addition morphism
$X^{(d_1)}\times X^{(d_2)}\to X^{(d)}$. For $\S_1\in D(X^{(d_1)})$, $\S_2\in D(X^{(d_2)})$,
let $\S_1\star \S_2\in D(X^{(d)})$ denote the object $(sum_{d_1,d_2})_!(\S_1\boxtimes \S_2)$.

Let $E^{\boxtimes d}\in D(X^d)$ denote the external
power of $E$; this sheaf is naturally $\Sigma_d$-equivariant. Consider
$$E\star...\star E\star...\star E\simeq (\sym_d)_!(E^{\boxtimes d})\in D(X^{(d)}).$$
Since the map $\sym_d$ is $\Sigma_d$-invariant, this sheaf is $\Sigma_d$-equivariant.

Recall that if $\Sigma$ is a finite group (acting trivially on a variety $\Y$),
we have the derived functor of invariants
$$R\on{Inv}_\Sigma:D^+(\Y)^\Sigma\to D^+(\Y)^.$$
By applying this functor to $(\sym_d)_!(E^{\boxtimes d})$ we obtain an object
$R\on{Inv}_{\Sigma_d}\left((\sym_d)_!(E^{\boxtimes d})\right)\in D^+(X^{(d)})$.

Finally, we set $E^{(d)}$ to be the $0$-the cohomology {\it in the usual t-structure}
of the complex $R\on{Inv}_{\Sigma_d}\left((\sym_d)_!(E^{\boxtimes d})\right)$. In other words,
$E^{(d)}$ is obtained from $(\sym_d)_!(E^{\boxtimes d})$ by taking non-derived 
$\Sigma_d$-invariants in {\it the abelian category of sheaves}.

Since the functor of stalks is exact on the category of sheaves, we obtain that
for a point $\imath_D:\on{pt}\to X^{(d)}$, where $D=\Sigma\, d_i\cdot x_i$ 
is an effective divisor of degree $d$ with the $x_i$'s pairwise distinct,
$$\imath_D^*(E^{(d)})\simeq \underset{i}\otimes\, Sym^{!,d_i}(E_{x_i}),$$
where $E_{x_i}$ denotes the stalk of $E$ at $x_i$.

For two positive integers $d_1,d_2$ recall the subset
$(X^{(d_1)}\times X^{(d_2)})_{disj}\subset X^{(d_1)}\times X^{(d_2)}$.
We have:
\begin{equation}  \label{symmetric factorization}
sum_{d_1,d_2}^*(E^{(d)})|_{(X^{(d_1)}\times X^{(d_2)})_{disj}}\simeq
(E_1^{(d_1)}\boxtimes E_2^{(d_2)})|_{(X^{(d_1)}\times X^{(d_2)})_{disj}}.
\end{equation}

\medskip

Let $j:\overset{\circ}X{}^{(d)}\to X^{(d)}$ denote the embedding of the complement to the diagonal
divisor. It is easy to see that $\overset{\circ}E{}^{(d)}:=E^{(d)}|_{\overset{\circ}X{}^{(d)}}$
is a local system. 

\begin{prop}  \label{symmetric power is perverse}
The complex $E^{(d)}[d]$ is a perverse sheaf; moreover,
$$E^{(d)}[d]\simeq j_{!*}(\overset{\circ}E{}^{(d)}[d]).$$
\end{prop}

Note that the proposition implies that the construction of $E^{(d)}$ is essentially
Verdier self-dual, i.e., ${\mathbb D}(E^{(d)}[d])\simeq (E^{*})^{(d)}[d]$,
where $E^{*}$ is the dual local system. This is so, because the isomorphism
obviously holds over $\overset{\circ}X{}^{(d)}$, from which both sides are extended
minimally. In particular, we have a canonical projection $(\sym_d)_!(E^{\boxtimes d})\to
E^{(d)}$.

The above fact about self-duality implies the following description of the 
co-stalks of $E^{(d)}$: For $D=\Sigma\, d_i\cdot x_i$ as above,
$$\imath_D^!(E^{(d)})\simeq \underset{i}\otimes\, Sym^{*,d_i}(E_{x_i})[-2d].$$

\medskip

Note also that from the proposition it follows that $E^{(d)}[d]$ embeds into
the perverse sheaf
$h^0_{perv}\left(R\on{Inv}_{\Sigma_d}\left((\sym_d)_!(E^{\boxtimes d})\right)[d]\right)$,
where $h^0_{perv}$ denotes the functor of taking $0$-th cohomology in the perverse
t-structure. However, this map is not in general an isomorphism. 

Indeed, let $\bF$
be of characteristic $2$, take $d=2$ and $E$ to be the trivial $1$-dimensional local system.
Then, $(\sym_2)_!(E^{\boxtimes 2})\simeq (\sym_d)_!(\bF_{X^2})$, and we have a canonical
embedding $\bF_{X^{(2)}}\to (\sym_d)_!(\bF_{X^2})$, 
whose cone is $j_!(\bF_{\overset{\circ}X{}^{(2)}})$.
By duality, we have an embedding of perverse sheaves
$j_*(\bF_{\overset{\circ}X{}^{(2)}})[2]\to (\sym_d)_!(\bF_{X^2})[2]$,
and it is easy to see that $j_*(\bF_{\overset{\circ}X{}^{(2)}})[2]$ identifies with the subobject
of $\Sigma_2$-invariants in $(\sym_d)_!(\bF_{X^2})[2]$.

\medskip

To prove \propref{symmetric power is perverse} we will need the following lemma,
which follows from \lemref{linear algebra powers}:

\begin{lem}  \label{symmetric power direct sum}
For $E=E_1\oplus E_2$ we have a canonical isomorphism:
$$E^{(d)}\simeq \underset{d_1+d_2=d}\oplus\, E_1^{(d_1)}\star E_2^{(d_2)}.$$
\end{lem}

\begin{proof} (of \propref{symmetric power is perverse})

The question being \'etale-local, we can assume that the local system
$E$ is trivial. Hence we can decompose $E=E_1\oplus...\oplus  E_m$, where
$E_i$'s are $1$-dimensional. By \lemref{symmetric power direct sum}, we have:
$$E^{(d)}\simeq \underset{\db}\oplus\, E_1^{(d_1)}\star...\star E_m^{(d_m)},$$
where $\db=(d_1,...,d_m)$ runs over the set if $m$-tuples of non-negative integers
with $\underset{i}\Sigma\, d_i=d$.

It is easy to see that for the trivial $1$-dimensional local system 
its $d$-th symmetric power is the trivial $1$-dimensional local system on
$X^{(d)}$. This implies that assertion of the proposition, since for each $\db$,
the map $X^{(d_1)}\times...X^{(d_m)} \to X^{(d)}$ is finite.

\end{proof}

\ssec{}   \label{sect ext ext}

We will now construct another complex of sheaves, denoted $\Lambda^{!,(d)}(E)$, on
$X^{(d)}$, called the external exterior power of $E$.

\begin{thm}    \label{external exterior}
For every $d\geq 1$ there exists a canonically defined complex 
$\Lambda^{!,(d)}(E)\in D(X^{(d)})$ endowed with a map
$\alpha_d:\Lambda^{!,(d)}(E)\to \Lambda^{!,(d-1)}(E)\star E$, such that:

\smallskip

\noindent{\em (1)}
The restriction $\overset{\circ}\Lambda{}^{!,(d)}(E):=
j^*(\Lambda^{!,(d)}(E))$ is a local system and 
$$\overset{\circ}\Lambda{}^{!,(d)}(E) \simeq
R\on{Inv}_{\Sigma_d}\left(j^*\left((\on{sym}_d)_!(E^{\boxtimes d})\right)\otimes sign\right).$$

\smallskip

\noindent{\em (2)}
For $d_1+d_2=d$ the restriction of
$sum_{d_1,d_2}^*(\Lambda^{!,(d)}(E))$ to 
$(X^{(d_1)}\times X^{(d_2)})_{disj}$ is canonically isomorphic
to the restriction to this open subset of
$\Lambda^{!,(d_1)}(E)\boxtimes \Lambda^{!,(d_2)}(E)$.

\smallskip

\noindent{\em (3)}
$\Lambda^{!,(d)}(E)[d]$ is perverse.

\smallskip

\noindent{\em (4)}
The composition
$$\Lambda^{!,(d)}(E)\overset{\alpha_d}\longrightarrow \Lambda^{!,(d-1)}(E)\star E
\overset{\alpha_{d-1}}\longrightarrow \Lambda^{!,(d-2)}(E)\star E\star E\to
\Lambda^{!,(d-2)}(E)\star E^{(2)}$$
is zero and the resulting complex of perverse sheaves 
\begin{align*}
&\Lambda^{!,(d)}(E)[d]\to \Lambda^{!,(d-1)}(E)\star E[d]\to...\to
\Lambda^{!,(i)}(E)\star E^{(d-i)}[d]\to \\
&\Lambda^{!,(i-1)}(E)\star E^{(d-i+1)}[d]\to...\to
E\star E^{(d-1)}[d]\to E^{(d)}[d]
\end{align*}
is exact.

\smallskip

\noindent{\em (5)}
For a divisor $D=\Sigma\, d_i\cdot x_i$ with the $x_i$'s pairwise distinct,
the co-stalk $\imath_D^!(\Lambda^{!,(d)}(E))$ is (quasi-) isomorphic to
$\underset{i}\otimes\, \Lambda^{!,d_i}(E_{x_i})[-2d]$, so that the 
co-stalk of the complex of point (4) identifies with the product over $i$
of the Koszul complexes of \lemref{linear algebra powers}(3).

\end{thm}

Note that the construction of $\Lambda^{!,(d)}(E)$ is not Verdier self-dual.
Set $$\Lambda^{*,(d)}(E):={\mathbb D}\left(\Lambda^{!,(d)}(E^{*})\right)[-2d].$$
Then $\Lambda^{*,(d)}(E)$ would satisfy the same properties (1),(2) and (3) of
\thmref{external exterior} as $\Lambda^{!,(d)}(E)$. Instead of point (4) we will
have an exact complex
\begin{align*}
&E^{(d)}[d]\to E^{(d-1)}\star E[d]\to...\to
E^{(i)}\star \Lambda^{*,(d-i)}(E)[d] \to \\
&E^{(i-1)}\star \Lambda^{*,(d-i+1)}(E)[d] \to...\to
E\star \Lambda^{*,(d-1)}(E)[d]\to \Lambda^{*,(d)}(E)[d].
\end{align*}
Instead of point (5), we would be able to describe the stalks of
$\Lambda^{*,(d)}(E)$:
$$\imath_D^*(\Lambda^{*,(d)}(E))\simeq 
\underset{i}\otimes\, \Lambda^{*,d_i}(E_{x_i}).$$

Observe also that when $\on{char}(\bF)=0$, both $\Lambda^{!,(d)}(E)$ and $\Lambda^{*,(d)}(E)$
are isomorphic to the minimal extension $j_{!*}(\overset{\circ}\Lambda{}^{!,(d)}(E))$.

\ssec{Proof of \thmref{external exterior}}

We proceed by induction on $d$. Evidently, for $d=1$ we can take 
$\Lambda^{!,(d)}(E)=E\in D(X)$. Thus, we can assume that $\Lambda^{!,(i)}(E)$
satisfying conditions (1)-(5) of \thmref{external exterior} have been constructed
for $i<d$.

Define $\Lambda^{!,(d)}(E)\in D(X^{(d)})[d]$ to be represented by the complex
{\it of perverse sheaves}
$$\K(d,E):=\Lambda^{!,(d-1)}(E)\star E[d]\to...\to
\Lambda^{!,(i)}(E)\star E^{(d-i)}[d]\to...\to E\star E^{(d-1)}[d]\to E^{(d)}[d].$$
(Here we are using the fact that the category of complexes of perverse sheaves
maps to $D(X^{(d)})$;  in fact, due to a theorem of Beilinson, 
the corresponding functor from the derived category of perverse sheaves
to $D(X^{(d)})$ is an equivalence.)

It is easy to see that $\Lambda^{!,(d)}(E)$ satisfies conditions (1) and (2).
Let us show that $\Lambda^{!,(d)}(E)[d]$ is a perverse sheaf. Since the question
is \'etale-local, we can assume that $E$ is the trivial local system. and let us
write $E=E_1\oplus...\oplus E_m$, where the $E_i$'s are trivial $1$-dimensional.

\begin{lem}
For $E=E_1\oplus E_2$ we have a canonical isomorphism:
$$\Lambda^{!,(d)}(E)\simeq \underset{d_1+d_2=d}\oplus\, \Lambda^{!,(d_1)}(E_1)\star 
\Lambda^{!,(d_2)}(E_2).$$
\end{lem}

\begin{proof}
Assume the validity of the lemma for $d'<d$. Then, by the induction hypothesis
and \lemref{symmetric power direct sum}, we obtain an isomorphism of complexes
of perverse sheaves:
$$\K(d,E)\simeq \underset{d_1+d_2=d}\oplus\, \K(d_1,E_1)\star \K(d_2,E_2).$$
This implies our assertion.
\end{proof}

Thus, by decomposing $E$ as a direct sum $E=E_1\oplus...\oplus E_m$ with the
$E_i$'s being $1$-dimensional, we reduce the perversity assertion to the case when
$E$ is itself $1$-dimensional. In the latter case we claim that
$$\Lambda^{!,(d)}(E)\simeq j_*(\overset{\circ}\Lambda{}^{!,(d)}(E)).$$
Indeed, by point (5) and the induction hypothesis, the co-stalk of
$\Lambda^{!,(d)}(E)$ at $D=\Sigma\, d_i\cdot x_i$ is quasi-isomorphic to
$$\underset{i}\otimes\, \Lambda^{!,d_i}(E_{x_i})[-2d].$$
Therefore, if one of the $d_i$'s is $>1$, then the corresponding
$\Lambda^{!,d_i}(E_{x_i})=0$, and the above expression vanishes.
Moreover, this shows that the constructed complex $\Lambda^{!,(d)}(E)$ 
satisfies condition (5) of the theorem.

Hence, the perversity assertion follows from the fact that the embedding
$j:\overset{\circ}X{}^{(d)}\to X^{(d)}$ is affine.

\medskip

Note that the last point of the proof shows $\Lambda^{!,(d)}(E)$ is not a sheaf
in the usual t-structure (cf. example preceding the proof of 
\propref{symmetric power is perverse}). However, the object $\Lambda^{*,(d)}(E)$
is always a sheaf.

\section{Appendix A: Proof of \thmref{strong form}}

\ssec{}
The assertion of \thmref{strong form} divides into two parts.
First, we will show that \thmref{main} remains valid under the assumption
that $\on{char}(\bF)\neq 2$. Secondly, we will prove a particular case of
\conjref{GL_n Hecke} under the assumption that $\S_E$ is a cuspidal
perverse sheaf.

The ingredient in the proof of \thmref{main} that relied on the assumption
that $\on{char}(\bF)> 2n$ was the proof of \thmref{vanishing}. The latter used
this assumption in the following two places:

\medskip

\noindent 1) Proof of the fact that if $E$ is irreducible of rank $>n$, then the
functor $\on{Av}^1_E$ is exact on $\wt{D}(\Bun_n)$.

\smallskip

\noindent 2) Proof of \lemref{aux hecke}.

\medskip

Let us first treat point 2). Let us call a representation $V\in \on{Rep}(\check G_\bF)$
positive if the action of the group $GL_n$ on it extends to an action of
the semi-group $\on{Mat}_{n,n}$. Every positive $V$ can be decomposed into
a direct sum of $V_d$, $d\geq 0$, according to the action of the center.
Equivalently, $V$ is positive of degree $d$ if it can be realized as a quotient
of (several copies of) $V_0^{\otimes d}$.

Let $\Gr^+_{GL_n,x}\subset \Gr_{G,x}$ be the "positive" part of the affine
Grassmannian, corresponding to the condition that the modification of vector
bundles $\beta:\M\to \M^0$ is such that $\M^0\to \M$ is a regular map of
coherent sheaves. One can show that $V$ is positive if and only if the
corresponding perverse sheaf on $\Gr_{G,x}$ is supported on $\Gr^+_{GL_n,x}$.

Clearly, if $V$ is any finite-dimensional representation of $GL_n$, by
multiplying it by a sufficiently high power of the determinant, we can make
it positive. 

Therefore, to prove \lemref{aux hecke} it suffices to show that the functor
$H(V,\cdot)$ is well-defined and exact on $\wt{D}(\Bun_n)$ for $V$ positive.
However, the proof of the latter fact can be obtained by the same argument 
as the proof of this statement for $V=V_0$ in \cite{Ga}, Sect. 7.

\medskip

\ssec{}   
Let us now treat point 1) above. The assumption on $\on{char}(\bF)$ was used
in the proof of Theorem 2.14 of \cite{Ga} in the following situation:

Recall that if $\Sigma$ is a finite group, we have the abelian category
$\wt{P}{}^\Sigma(\Bun_n\times X)$, which is a Serre quotient of the
category of $\Sigma$-equivariant perverse sheaves on $\Bun_n\times X$,
with the group $\Sigma$ acting trivially. 

We will take $\Sigma=\Sigma_i$, and consider the functor
$$\S\mapsto (\S\otimes sign)_{\Sigma_i}$$
of $\Sigma_i$-anti-coinvariants. 

This functor is right-exact on
$P^{\Sigma_i}(\Bun_n\times X)$, and hence also on $\wt{P}{}^\Sigma(\Bun_n\times X)$,
and we need to replace by it the (exact, because of the assumption on
the characteristic) functor $\S\mapsto \on{Hom}_{\Sigma_i}(sign,\S)$ considered
in Sect. 3 of \cite{Ga}.

To make the argument work, we need to insure that an analog of Proposition 1.11 of 
\cite{Ga} holds in our situation. Namely, let us consider
the functor $\wt{P}(\Bun_n\times S)\to \wt{P}{}^{\Sigma_i}(\Bun_n\times X\times S)$
given by
$$\S\mapsto H^{\boxtimes i}_S(\S)|_{\Bun_n\times \Delta(X)\times S}[1-i].$$
This functor maps to the abelian category because of Property 1 of 
$\wt{D}(\Bun_n)$ and the ULA assertion (cf. Lemma 3.7 of \cite{Ga}). We have:

\begin{prop}
If $\on{char}(\bF)\neq 2$, then the right-exact functor
$$\S\mapsto \left(\left(H^{\boxtimes i}_S(\S)|_{\Bun_n\times \Delta(X)\times S}[1-i]\right)
\otimes sign\right)_{\Sigma_i}$$
mapping $\wt{P}(\Bun_n\times S)$ to $\wt{P}(\Bun_n\times X\times S)$ is $0$, if $i>n$,
and is isomorphic to $\S\mapsto m^*(\S)[1]$ if $i=n$.
\end{prop}

\begin{proof}

Using \lemref{aux hecke}, the functor considered in the proposition is isomorphic
to $$\S\mapsto H((V_0\otimes sign)^{\otimes i}_{\Sigma_i},\S).$$ Hence, the assertion
follows the fact that if $\on{char}(\bF)\neq 2$, then 
$(V_0\otimes sign)^{\otimes i}_{\Sigma_i}\simeq \Lambda^i(V_0)$.

\end{proof}

\ssec{}
Finally, let us show that if $\S_E$ is a cuspidal perverse sheaf, which
has a $GL_n$-Hecke property with respect to a local system $E$, then, in fact,
is satisfies the full Hecke property.

First, note that it suffices to construct the functorial isomorphisms
$\alpha(V)$ of \eqref{eigencondition}, and verify their properties, 
for representations $V$, which are positive.

We begin with the following observation:

\begin{lem}  \label{A.5}
For any $V\in \on{Rep}(\check G_\bF)$, and $\S\in D(\Bun_n)$, which is cuspidal
and perverse, the object $H(V,\S)$ is a perverse sheaf.
\end{lem}

\begin{proof}

Since the assertion is essentially Verdier self-dual, it suffices to show
that $H(V,\S)$ belongs to $D^{\leq 0}(\Bun_n\times X)$. 

Suppose the contrary, and consider the truncation map
\begin{equation} \label{trunk}
H(V,\S)\to \tau^{>0}(H(V,\S)).
\end{equation}
By \lemref{aux hecke}, $H(V,\S)$ is exact on $\wt{D}(\Bun_n\times X)$.
By assumption, $\S$ is cuspidal, which implies that $H(V,\S)$ is also
cuspidal.

Property 2 of $\wt{D}(\Bun_n\times X)$ (cf. \cite{Ga}, Sect. 2.12)
implies that the truncation map
\eqref{trunk} is $0$, which is a contradiction.

\end{proof}

Let $\underline{V^*_0}$ denote the vector space underlying the corresponding
representation of $\check G_\bF=GL_n$. Consider the object
$$H\left((V\otimes \underline{V^*_0})\boxtimes...\boxtimes 
(V\otimes \underline{V^*_0}),\S_E\right)
\in \on{D}(\Bun_n\times X^d).$$
It is $\Sigma_d$-equivariant, and by assumption, it is isomorphic to the perverse sheaf
$$\S_E\boxtimes \left((E\otimes \underline{V^*_0}[1])\boxtimes...
\boxtimes (V\otimes \underline{V^*_0}[1])\right).$$

By restricting both sides to the diagonal $\Bun_n\times X\subset \Bun_n\times X^d$,
we obtain a $\Sigma_d$-equivariant isomorphism
\begin{equation}  \label{isom of tensor}
H\left((V\otimes \underline{V^*_0})^{\otimes d},\S_E\right)\simeq \S_E\boxtimes
(E\otimes \underline{V^*_0})^{\otimes d}[1].
\end{equation}
Moreover, both sides of \eqref{isom of tensor} are acted on by the group $GL(\underline{V_0})$
of automorphisms of the vector space $\underline{V_0}$, and the isomorphism of
\eqref{isom of tensor} is compatible with these actions.

Let us take $\Sigma_d$-coinvariants of both sides of \eqref{isom of tensor}.
By \lemref{A.5}, we obtain a $GL(\underline{V_0})$-equivariant isomorphism of
perverse sheaves:
\begin{equation} \label{isom of sym}
H\left(\Sym^{*,d}(V\otimes \underline{V^*_0}),\S_E\right)\simeq
\S_E\boxtimes \Sym^{*,d}(E\otimes \underline{V^*_0})[1].
\end{equation}

Let now $V\in \on{Rep}(\check G_\bF)$ be a positive representation of degree $d$.
Let $\underline{V}$ denote the underlining vector space, which we may regard
as a representation of $GL(\underline{V_0})$. Note that we have an isomorphism
of $\check G_\bF$-representations
$$\left(\Sym^{*,d}(V\otimes \underline{V^*_0})\otimes \underline{V}\right)^{GL(\underline{V_0})}
\simeq V.$$
Let us tensor both sides of \eqref{isom of sym} by $\underline{V}$ and take
$GL(\underline{V_0})$-invariant parts. By \lemref{A.5}, we obtain:
$$H(V,\S_E)\simeq \S_E\boxtimes \left(\Sym^{*,d}(E\otimes \underline{V^*_0})\otimes
\underline{V}\right)^{GL(\underline{V_0})}[1]\simeq \S_E\boxtimes E^V[1],$$
where $E^V$ is the local system corresponding to $E$ and the $\check G_\bF$-representation
$V$.

Thus, we have constructed a functorial isomorphism $\alpha(V)$ for $V$
positive of a fixed degree. To check the commutativity of the diagrams
1) and 2) it is sufficient to do this when
$V_1\simeq (V\otimes \underline{V^*_0})^{d_1}$
and $V_2\simeq (V\otimes \underline{V^*_0})^{d_2}$. In this case, the required
commutativity follows by construction.

\section{Appendix B: Proof of \thmref{extended satake}}

In this section we will show how to deduce \thmref{extended satake} from
\thmref{satake}. We will work with sheaves over any ring of coefficients (e.g.,
$\bF$ or $\bO$) and we will regard the Langlands dual group $\cG$
as a group-scheme over this ring. We shall denote by $\Rep(\cG)$ the category
of algebraic representations of $\cG$.

\ssec{}

We shall first consider the case $d=1$. By \cite{MV}, Proposition 2.2, to any object $V\in \Rep(\cG)$
we can attach a spherical perverse sheaf $\CT_V$ on $\Gr_{G,X}$. Let $R$ be the algebra
of functions on $G$, viewed as an ind-object of $\Rep(\cG)$ via the left action of $\cG$
on itself. Let us denote by $\CR_X$ the corresponding ind-object of $P^{\CG}(\Gr_{G,X})$. 
The right action of $\cG$ on itself endows $\CR_X$ with a $\cG$-action.

We define the functor
$\sF:P^{\cG,1}\to \Ind(P^{\CG}(\Gr_{G,X}))$ by the formula
$$\CK\mapsto \left(s^*(\CK)\otimes \CR_X[1]\right)^{\cG},$$
where the superscript $\cG$ designates $\cG$-invariants. We shall
denote by the same symbol the extension of this functor onto $\Ind(P^{\cG,1})$.

\medskip

Let ${\bf 1}_{\Gr_{G,X}}$ denote the natural section $X\mapsto \Gr_{G,X}$. We define
the functor $\sG:P^{\CG}(\Gr_{G,X})\to \Ind(P^{\cG,1})$ by
$$\CT\mapsto h^0\left({\bf 1}_{\Gr_{G,X}}^!(\CR_X\star \CT)\right),$$
where $\star$ is the convolution functor on $P^{\CG}(\Gr_{G,X})$, and $h^0\left(\cdot\right)$
designates the $0$-th perverse cohomology. We will denote by the same symbol
the extension of $\sG$ onto $\Ind(P^{\CG}(\Gr_{G,X}))$.

\begin{prop}  \label{d=1} 
The functors $\sF$ and $\sG$ map $P^{\cG,1}$ to $P^{\CG}(\Gr_{G,X})$
and $P^{\CG}(\Gr_{G,X})$ to $P^{\cG,1}$, respectively, and define mutually inverse
equivalences of categories.
\end{prop}

The rest of this subsection is devoted to the proof of this proposition. 

\medskip

Let us consider the composition $\sG\circ \sF:P^{\cG,1}\to \Ind(P^{\cG,1})$. Since the
functor of convolution is exact, for $\CK\in P^{\cG,1}$ we have:
$$\CR_X\star \left(s^*(\CK)\otimes \CR_X[1]\right)^{\cG}\simeq
\left(s^*(\CK)[1]\otimes (\CR_X\star \CR_X)\right)^\cG,$$
where $\cG$ acts on $\CR_X\star \CR_X$ via the second multiple.

However, $\CR_X\star \CR_X\simeq \CR_X\otimes R$, with the diagonal
action of $\cG$. Hence, we must calculate
\begin{equation} \label{prelcomp 1}
h^0\left({\bf 1}_{\Gr_{G,X}}^!(s^*(\CK\otimes R)[1]\otimes\CR_X)\right).
\end{equation}

\begin{lem}   \label{pre comp 1}
For any $\CK'\in \Ind(P^{\cG,1})$ we have a canonical isomorphism
$$\CK'\simeq h^0\left({\bf 1}_{\Gr_{G,X}}^!(s^*(\CK')[1]\otimes\CR_X)\right)$$
\end{lem}

\begin{proof}

The embedding of the trivial representation into $R$ gives rise to a map
\begin{equation} \label{delta}
{\bf 1}_{\Gr_{G,X}}{}_*(\on{Const}_X[1])\to \CR_X,
\end{equation}
where $\on{Const}_X$ denotes the constant sheaf on $X$. Hence, for any $\CK'$
as above, we have a map
$$\CK'\to {\bf 1}_{\Gr_{G,X}}^!(s^*(\CK')[1]\otimes\CR_X),$$
and we must show that the LHS identifies with the maximal sub-object 
of the RHS, supported on the image of ${\bf 1}_{\Gr_{G,X}}$.

Let us first assume that $\CK'$ is lisse. Then $s^*(\CK')\otimes \CR_X\in
\Ind(P^{\CG}(\Gr_{G,X}))$ is ULA with respect to the projection $s:\Gr_{G,X}\to X$.
Therefore, to prove our assertion,  it would be sufficient to show that for some 
(or any) point $x\in X$,
\begin{equation} \label{fiber of functor}
\CK'_x\to h^0\left({\bf 1}_{\Gr_{G,x}}^!(\CK'_x\otimes\CR)\right)
\end{equation}
is an isomorphism, where $\CK'_x$ denotes the fiber of $\CK'$ at $x$,
and ${\bf 1}_{\Gr_{G,x}}^!$ is the embedding of the unit point into $\Gr_{G,x}$.
However, the latter assertion follows from the equivalence
$\Rep(\CG)\simeq P^{G(\CO_x)}(\Gr_{G,x})$. 

\medskip

Let now $\CK'$ be arbitrary. We can write
\begin{equation} \label{decomp}
\CK'_1\to \CK'\to j_*j^*(\CK')\to \CK'_2,
\end{equation}
where $j$ is the embedding of an open subset $X'\hookrightarrow X$,
and $\CK'_1,\CK'_2$ are perverse sheaves supported on $X-X'$.
By choosing $X'$ to be sufficiently small, we can arrange that the
restriction of $\CK'$ to $X'$ be lisse.

Since the functor $h^0\left({\bf 1}_{\Gr_{G,x}}^!(\cdot)\right)$
is left exact, it is sufficient to show that the map in the lemma
is an isomorphism for $\CK'_1$, $\CK'_2$ and $j_*j^*(\CK')$. 
The assertion concerning $\CK'_1$ and $\CK'_2$ is an immediate
corollary of the equivalence $\Rep(\cG)\simeq P^{G(\CO_x)}(\Gr_{G,x})$.
In addition, we have:
$${\bf 1}_{\Gr_{G,X}}^!(s^*(j_*j^*(\CK'))\otimes\CR_X)\simeq
j_*j^*\left({\bf 1}_{\Gr_{G,X}}^!(s^*(\CK')\otimes\CR_X)\right),$$
and our assertion follows from the lisse case, considered above.

\end{proof}

Using the lemma, we obtain 
$$\sG\circ \sF(\CK)\simeq \left(\CK\otimes R\right)^\cG\simeq \CK,$$
as required. Let us now consider the composition $\sF\circ \sG:P^{\CG}(\Gr_{G,X})\to
\Ind(P^{\CG}(\Gr_{G,X}))$.

\begin{lem} \label{pre comp 2}
For $\CT\in P^{\CG}(\Gr_{G,X})$ there exists a canonical isomorphism
$$s^*\left(h^0\left({\bf 1}_{\Gr_{G,X}}^!(\CR_X\star \CT)\right)[1]\right)\otimes \CR_X
\simeq \CR_X\star \CT.$$
\end{lem}

\begin{proof}

Let us rewrite the LHS of the expression in the lemma as
\begin{equation} \label{exp 2}
\CR_X\star {\bf 1}_{\Gr_{G,X}}{}_*
\left(h^0\left({\bf 1}_{\Gr_{G,X}}^!(\CR_X\star \CT)\right)\right).
\end{equation}

The natural map
$${\bf 1}_{\Gr_{G,X}}{}_*
\left(h^0\left({\bf 1}_{\Gr_{G,X}}^!(\CR_X\star \CT)\right)\right)\to \CR_X\star \CT$$
gives rise to a map from the expression in \eqref{exp 2} to 
$$\CR_X\star (\CR_X\star \CT)\simeq (\CR_X\star \CR_X)\star \CT.$$
The multiplication on $R$ give rise to a map $\CR_X\star \CR_X\to \CR_X$.
Hence, by composing, we obtain a map 
$$\CR_X\star {\bf 1}_{\Gr_{G,X}}{}_*
\left(h^0\left({\bf 1}_{\Gr_{G,X}}^!(\CR_X\star \CT)\right)\right)\to \CR_X\star \CT.$$

To show that the latter map is an isomorphism, we proceed as in the proof of
\lemref{pre comp 1}, by reducing the assertion to the case when $\CT$ is ULA
with respect to $s:\Gr_{G,X}\to X$.

\end{proof}

Thus,
$$\sF\circ \sG(\CT)\simeq \left(\CR_X\star \CT\right)^{\cG},$$
where $\cG$ acts on the convolution via its action on $\CR_X$. The map
\eqref{delta} gives rise to a map
$$\CT\to \left(\CR_X\star \CT\right)^{\cG},$$
and repeating the argument used in the proofs of Lemmas \ref{pre comp 1} and \ref{pre comp 2},
we show that the latter map is an isomorphism.

Thus, we obtain that the functors $\sF$ and $\sG$ induce mutually inverse equivalences
of categories of $\Ind(P^{\CG}(\Gr_{G,X}))$ and $\Ind(P^{\cG,1})$. In particular, 
$\sF$ sends $P^{\cG,1}$ to $P^{\CG}(\Gr_{G,X})$ and $\sG$ sends $P^{\CG}(\Gr_{G,X})$
to $P^{\cG,1}$, and the resulting functors $P^{\cG,1}\rightleftarrows P^{\CG}(\Gr_{G,X})$
are also mutually inverse.

\ssec{}

To treat the case of an arbitrary $d$, we will have to construct an object
$^f\CR^{(d)}$ in $P^{\CG^d}(\Gr^d_G)$, which will play a role similar to that of
$\CR_X$. This construction is essentially borrowed from 
\cite{CHA}, Sect. 3.4.

First, let us recall the following construction from \cite{MV}, Sect. 5. Let
$\Delta$ denote the embedding of the diagonal $X\to X\times X$ and
$j$ be the embedding of the complement. Then given two objects
$V,W\in \Rep(\cG)$ there exists a canonical map
$$j_*j^*(\CT_V\boxtimes \CT_W)\to \Delta_*(\CT_{V\otimes W}).$$

Applying this to $V=W=R$, and composing with the
map $\CT_{R\otimes R}\to \CT_R$, corresponding to the algebra
structure on $R$, we obtain a map
\begin{equation}   \label{chiral}
j_*j^*(\CR_X\boxtimes \CR_X)\to \Delta_*(\CR_X).
\end{equation}

For a positive integer $d$, let us denote by $J^d$ the finite set 
$\{1,...,d\}$, and for a surjection of finite sets 
$J^d\twoheadrightarrow I$, let $\Delta^I$ be the embedding of
the corresponding diagonal $X^I\hookrightarrow X^d=X^{J^d}$. Let
$\oDelta^I$ be the locally closed embedding of the open subset
$\oX^I\hookrightarrow X^{J^d}$, obtained by removing from $X^I$
its diagonal divisor. Let us denote by $\Gr^I_G$ the corresponding
version of the affine Grassmannian over $X^I$ (obtained as a
pull-back from $\Gr^{|I|}_G$ over $X^{(|I|)}$, and let $\oGr^I_G$ be
its restriction to $\oX^I$. By a slight abuse of notation, we shall
denote by $\Delta^I$ and $\oDelta^I$ the embeddings
of $\Gr^I_G$ and $\oGr^I_G$ into $\Gr^{J^d}_G$, respectively.
The basic factorization property of the affine Grassmannian
(cf. \cite{MV}, Sect. 5) is that we have a canonical isomorphism
$$\oGr^I_G\simeq \left(\Gr_{G,X}\right)^{\times I}|_{\oX^I}.$$

\medskip

Consider the perverse sheaf $\CR_X^{\boxtimes I}|_{\oGr^I_G}$, and
let us denote by $\oR_X^I$ its *-direct image under
$\oDelta^I:\oGr^I_G\to \Gr_G^{J^d}$. From \eqref{chiral} we obtain that
for any surjection $I\twoheadrightarrow I'$ with $|I'|=|I|-1$ there exists
a naturally defined map
$$\oR_X^I\to \oR_X^{I'}.$$

By appropriately choosing the signs (cf. \cite{CHA}, 3.4.11), we obtain
a complex of perverse sheaves on $\Gr^{J^d}_G$, which we will denote by 
$\fC^\bullet(\CR^d_X)$, and whose $k$-th term is
$$\underset{I,|I|=d-k}\oplus\, \oR_X^I.$$

By essentially repeating the proof of Lemma 2.4.12 of \cite{CHA}, we obtain
the following:
\begin{lem}  \label{cousin}
The complex $\fC^\bullet(\CR^d_X)$ is acyclic off degree $0$.
\end{lem}

Let us denote by $^f\CR^d_X$ the $0$-th cohomology of $\fC^\bullet(\CR^d_X)$.

\medskip

\noindent {\it Remark.}
The perverse sheaf $^f\CR^d_X$ has been introduced by A.~Beilinson
in the construction of automorphic sheaves using a "spectral projector".

\ssec{}

We shall now construct a version of $^f\CR^d_X$ that lives on $\Gr_G^d$ instead
of $\Gr_G^{J^d}$, which we will denote by $^f\CR^{(d)}_X$. Let us denote by 
$\on{sym}_d$ the natural map $\Gr_G^{J^d}\to \Gr_G^d$. If we worked with 
sheaves with characteristic $0$ coefficients, we would define $^f\CR^{(d)}_X$
as $\Sigma_d$-invariants in ($\on{sym}_d)_!({}^f\CR^d_X)$. In the case of arbitrary
coefficients, we proceed as follows.

We define $^f\CR_X^{(d)}$ as the kernel of the map
$$\Bigl((\on{sym}_d)!\left(\fC^0(\CR^d_X)\right)\Bigr)_{\Sigma_d}\to
\Bigl((\on{sym}_d)!\left(\fC^1(\CR^d_X)\right)\Bigr)_{\Sigma_d},$$
where the subscript "$\Sigma_d$" means $\Sigma_d$-coinvariants,
taken in the category of perverse sheaves.

By construction, $^f\CR^{(d)}_X$ has the following factorization property.
For a partition $\ol{d}:d=d_1+...+d_m$, we have an isomorphism
\begin{equation} \label{factor factor}
^f\CR^{(d)}_X|_{X^{(\ol{d})}_{disj}\underset{X^{(d)}}\times \Gr_G^d}\simeq
\left(\CR^{(d_1)}_X\boxtimes...\boxtimes \CR^{(d_m)}_X\right)|
_{X^{(\ol{d})}_{disj}\underset{X^{(d)}}\times (\Gr_G^{d_1}\times...\times \Gr_G^{d_m})},
\end{equation}
under the identification
$$X^{(\ol{d})}_{disj}\underset{X^{(d)}}\times \Gr_G^d\simeq
X^{(\ol{d})}_{disj}\underset{X^{(d)}}\times (\Gr_G^{d_1}\times...\times \Gr_G^{d_m}).$$

\medskip

Consider the complex $\fC^\bullet(\CR^{(d)}_X)$, whose terms are given by
$$\fC^k(\CR^{(d)}_X):=\Bigl((\on{sym}_d)!\left(\fC^k(\CR^d_X)\right)\Bigr)_{\Sigma_d}.$$

Note that the $k$-term is a *-extension of a complex on $\Gr^d_G$, which is supported
over the locus of codimension $k$ in $X^{(d)}$, corresponding the collision pattern
of points.

\begin{prop}
The complex 
$$^f\CR^{(d)}_X\to \fC^0(\CR^{(d)}_X)\to \fC^1(\CR^{(d)}_X)\to...\to \fC^{d-1}(\CR^{(d)}_X)$$
is acyclic.
\end{prop}

\begin{proof}

The assertion is evidently true for $d=2$ and we proceed by induction. By \eqref{factor factor},
we can assume that the complex in question is acyclic off the preimage of main diagonal
$\Delta:X\hookrightarrow X^{(d)}$. 

Since the last arrow $\fC^{d-2}(\CR^{(d)}_X)\to \fC^{d-1}(\CR^{(d)}_X)$ is surjective,
the assertion of the proposition is equivalent to the fact that the map
$$h^{d-1}\left(\Delta^!({}^f\CR^{(d)}_X)\right)\to \CR_X,$$
resulting from the above complex, is an isomorphism.

\medskip

Consider the perverse sheaf $\left((\on{sym}_d)_!({}^f\CR^d_X)\right)_{\Sigma_d}$
on $\Gr^d_G$, which maps naturally to $^f\CR^{(d)}_X$. Both the kernel and the cokernel
of this map are supported over the preimage of the diagonal divisor in $X$.
Therefore, the map
\begin{equation} \label{map on diag}
h^{d-1}\Bigl(\Delta^!\left((\on{sym}_d)_!({}^f\CR^d_X)\right)_{\Sigma_d}\Bigr)\to
h^{d-1}\left(\Delta^!({}^f\CR^{(d)}_X)\right)
\end{equation}
is surjective. 

However, it is easy to see that 
$$h^{d-1}\Bigl(\Delta^!\left((\on{sym}_d)_!({}^f\CR^d_X)\right)_{\Sigma_d}\Bigr)\simeq 
\Delta_*(\CR_X),$$
and the resulting composition
$$\CR_X\to h^{d-1}\left(\Delta^!({}^f\CR^{(d)}_X)\right)\to \CR_X$$
is the identity map. Hence, the map of \eqref{map on diag} is an isomorphism.

\end{proof}

\ssec{}

Let us introduce the hybrid category $\Rep(\cG,P^{\CG^d}(\Gr^d_G))$,
which consists of $\CG^d$-equivariant ind-perverse sheaves $\CT$ on 
$\Gr^d_G$, such that for every partition $\ol{d}:d=d_1,...,d_n$
the pull-back of $\CT$ to
$X^{(\ol{d})}_{disj}\underset{X^{(d)}}\times (\Gr_G^{d_1}\times...\times \Gr_G^{d_m})$
carries an action of $\cG^{\times m}$, such that conditions, parallel to 1) and 2) in
the definition of $P^{\cG,d}$ hold. By construction, $^f\CR_X^{(d)}$ is an object of $\Rep(\cG,P^{\CG^d}(\Gr^d_G))$.

\medskip

We define the functor
$$\CT'\mapsto \CT'{}^\cG:\Rep(\cG,P^{\CG^d}(\Gr^d_G))\to
\Ind(P^{\CG^d}(\Gr^d_G)),$$ 
which sends $\CT'$ to its maximal sub-ind perverse sheaf, on which all
the actions of $\cG^{\times m}$ are trivial. 

For an object $\CK\in P^{\cG,d}$ consider 
$$s^*(\CK)\shriektimes {}^f\CR_X[-d]\in \Rep(\cG,P^{\CG^d}(\Gr^d_G)),$$
where $s$ denotes the projection $\Gr^d_G\to X^{(d)}$, and where
$\shriektimes$ denotes the functor ${\mathbb D}({\mathbb D}
(\cdot) \otimes {\mathbb D}(\cdot))$. 

We define the functor $\sF^d:P^{\cG,d}\to \Ind(P^{\CG^d}(\Gr^d_G))$ by
$$\CT\mapsto  \left(s^*(\CK)\shriektimes {}^f\CR_X[-d]\right)^\cG.$$

\medskip

Let ${\bf 1}_{\Gr^d_G}$ denote the unit section of $\Gr_G^d$. If $\CT'$ 
is an object of $\Rep(\cG,P^{\CG^d}(\Gr^d_G))$, then 
$h^0\left({\bf 1}_{\Gr^d_G}^!(\CT')\right)$ is naturally an object of $\Ind(P^{\cG,d})$.

For $\CT\in P^{\CG^d}(\Gr^d_G)$, we can consider $^f\CR_X^{(d)}\underset{\Gr^d_G}\star \CT$ 
as an object of $\Rep(\cG,P^{\CG^d}(\Gr^d_G))$, where $\underset{\Gr^d_G}\star$
refers to the convolution on $P^{\CG^d}(\Gr^d_G)$ (as opposed to the one, involving
$d=d_1+d_2$, as was the case in \secref{rev satake}).

We define the functor 
$\sG^d:P^{\CG^d}(\Gr^d_G)\to \Ind(P^{\cG,d})$ by
$$\CT\mapsto h^0\left({\bf 1}_{\Gr^d_G}^!({}^f\CR_X^{(d)}\star \CT)\right).$$

\begin{prop}
The functors $\sF^d$ and $\sG^d$ map to $P^{\CG^d}(\Gr^d_G)$ and
$P^{\cG,d}$, respectively, and are mutually inverse equivalences of categories.
\end{prop}

The proof of this proposition essentially repeats that of \propref{d=1}, using the
properties of $^f\CR_X^{(d)}$ established in the previous subsection. This establishes
the equivalence $P^{\cG,d}\simeq P^{\CG^d}(\Gr^d_G)$, stated in Theorem 2.6.
It remains to show the compatibility of this equivalence and the $\star$ operations
in (2) and (3).

For $d=d_1+d_2$, we have a natural morphism $^f\CR_X^{(d_1)}\star {}^f\CR_X^{(d_2)}\to
{}^f\CR_X^{(d)}$. It gives rise to a functorial morphism, defined for 
$\CK_1\in P^{\cG,d_1}$, $\CK_2\in P^{\cG,d_2}$:
$$\sF^{d_1}(\CK_1)\star \sF^{d_2}(\CK_2)\to \sF^d(\CK_1\star \CK_2).$$

To show that the latter morphism is an isomorphism, we decompose $\Gr^d_G$
with respect to the diagonal stratification of $X^{(d)}$, and the assertion follows
from the fact that the isomorphism of \propref{d=1} intertwines convolution of
$\CG$-equivariant perverse sheaves with tensor products of objects of
$P^{\cG,1}$.


\begin{thebibliography}{99}

\bibitem{BBD} A.~Beilinson, J.~Bernstein and P.~Deligne,
{\em Faisceaux pervers}, Ast\'erisque {\bf 100} (1982).

\bibitem{CHA} A.~Beilinson, V.~Drinfeld, {\em Chiral algebras},
AMS Colloquium Publ. {\bf 51}, Providence (2004).

\bibitem{BM} W.~Borho, R.~MacPherson,
{\em Small resolutions of nilpotent varieties}, 
Ast\'erisque {\bf 101-102} (1982), 23-74.

\bibitem{DJ}
A~J.~de Jong, {\em A conjecture on arithmetic
fundamental groups}, Israel Jour. of Math. {\bf 121} (2001), 61--64.


\bibitem{FGV}
E.~Frenkel, D.~Gaitsgory and K.~Vilonen,
{\em On the geometric Langlands conjecture}, J. Amer. Math. Soc. {\bf 15} (2002), 
367--417.

\bibitem{Ga}
D.~Gaitsgory, {\em On a vanishing conjecture appearing in the geometric Langlands
correspondence}, Annals of Mathematics, {\bf 160} (2004), 617-682.

\bibitem{Ga1}
D.~Gaitsgory {\em Construction of central elements in the affine Hecke
algebra via nearby cycles}, Inv. Math., {\bf 144} (2001), 253--280. 

\bibitem{Ga2}
D.~Gaitsgory {\em Automorphic sheaves and Eisenstein series},
PhD Thesis, Tel Aviv University (1997).

\bibitem{Ill}
L.~Illusie, {\em Th\'eorie de Brauer et Caract\'eristique 
d'Euler-Poincar\'e (d'apr\`es P.~Deligne)}, 
Ast\'erisque {\bf 82-83} (1981), 161--172.

\bibitem{La}
G.~Laumon {\em Faisceaux automorphes pour $GL_n$: la premi\`ere
construction de Drinfeld}, alg-geom/9511004.

\bibitem{MV}
I.~Mirkovi\'c and K.~Vilonen {\em Geometric Langlands duality and
representations of algebraic groups over commutative rings,}
math.RT/0401222.

\end{thebibliography}
\end{document}